\documentclass[11pt]{article}
\usepackage{}
\usepackage{mathrsfs}
\usepackage{amsfonts}

\usepackage{CJK}
\usepackage{amsfonts,amssymb,amsmath,mathrsfs}
\usepackage{color,latexsym,amsfonts}

 \newcommand{\qed}{\hfill\rule{2mm}{3mm}\vspace{4mm}}
 \def\bitemize{\begin{itemize}}\def\eitemize{\end{itemize}}

 \newtheorem{theorem}{Theorem}[section]
 \newtheorem{lemma}[theorem]{Lemma}
 \newtheorem{corollary}[theorem]{Corollary}
 \newtheorem{proposition}[theorem]{Proposition}
 
 \newtheorem{Definition}[theorem]{Definition}
 \newtheorem{remark}[theorem]{Remark}
 \newtheorem{condition}[theorem]{Condition}

 \def\blemma{\begin{lemma}\sl{}\def\elemma{\end{lemma}}}
 \def\btheorem{\begin{theorem}\sl{}\def\etheorem{\end{theorem}}}
 \def\bcorollary{\begin{corollary}\sl{}\def\ecorollary{\end{corollary}}}
 
 \def\bproposition{\begin{proposition}\sl{}\def\eproposition{\end{proposition}}}

 \def\beqlb{\begin{eqnarray}}\def\eeqlb{\end{eqnarray}}
 \def\beqnn{\begin{eqnarray*}}\def\eeqnn{\end{eqnarray*}}

 \def\mbb{\mathbb}\def\mbf{\mathbf}

 \def\<{\langle}\def\>{\rangle}

 \def\ar{&\!\!}\def\qqquad{\qquad\qquad}

 \def\eqref#1{{\rm(\ref{#1})}}

 \def\proof{\noindent{\it Proof.~}}\def\qed{\hfill$\Box$\medskip}

\newfam\msbmfam\font\tenmsbm=msbm10\textfont
\msbmfam=\tenmsbm\font\sevenmsbm=msbm7

\scriptfont\msbmfam=\sevenmsbm

\def\IR{\mathbb{R}}

\newcommand{\de}{\delta}

\newcommand{\la}{\lambda}
\newcommand{\ep}{\varepsilon}

\newcommand{\bR}{\mathbb{R}}

\newcommand{\goto}{\rightarrow}

\begin{document}

\noindent{(Draft: 2016/09/12)}

\bigskip

\centerline{\Large\bf A distribution-function-valued SPDE}
 \smallskip
\centerline{\Large\bf and its applications \footnote{Supported by grants from
NSFC 11301020, 11401012, 11602003 and  NSERC. }}

\bigskip\bigskip

\centerline{\ Li Wang, Xu Yang\footnote{Corresponding author
(xuyang@mail.bnu.edu.cn).}, Xiaowen Zhou }

\bigskip\bigskip

{\narrower{\narrower

\noindent{\bf Abstract.}  In this paper we further study  the
stochastic partial differential equation first proposed by Xiong
\cite{X13}. Under localized conditions on its coefficients,
we prove a comparison theorem on its solutions and show
that the solution is in fact distribution-function-valued. We also
establish  pathwise uniqueness of the solution. As applications we
obtain the well-posedness of  martingale problems for two classes of
measure-valued diffusions: interacting super-Brownian motions and
interacting Fleming-Viot processes. Properties of the two
superprocesses such as the existence of density fields and the
survival-extinction behaviors are also studied.

\medskip


\medskip


\par}\par}
\medskip

\section{Introduction}

\setcounter{equation}{0}

It is well known that the density process
$\{X_t(x):t>0,x\in\mbb{R}\}$ of a one-dimensional binary branching super-Brownian
motion solves the following non-linear stochastic partial
differential equation (SPDE):
 \beqlb\label{1.01}
\frac{\partial}{\partial t}X_t(x)=\frac12\Delta
X_t(x)+\sqrt{X_t(x)}\dot{W}_t(x), \qquad t>0,~x\in\mbb{R},
 \eeqlb
where $\Delta$ denotes the Laplacian operator and $\{\dot{W}_t(x): t
\geq 0, x\in\mbb{R}\}$ is the derivative of a space-time Gaussian
white noise. This SPDE was first derived and studied independently
by Konno and Shiga \cite{KS88} and Reimers \cite{Rei89}. The weak
uniqueness of solution to \eqref{1.01} follows from that of a
martingale problem for  super-Brownian motion. But the strong
(pathwise) uniqueness of nonnegative solution to \eqref{1.01}
remains open even though it has been studied by many authors. The
main difficulty comes from the unbounded drift coefficient and the
non-Lipschitz diffusion coefficient. Progresses have been made in
considering modified forms of SPDE \eqref{1.01}. When the random
field $\{W_t(x): t\geq 0, x\in \mbb{R}\}$ is colored in space and
white in time, the strong uniqueness of nonnegative solution to
\eqref{1.01}
 and to more general equations were studied in \cite{MPS06,RS13,N14}. When
$\{W_t(x): t\geq 0, x\in \mbb{R}\}$ is a space-time Gaussian white
noise and solutions are allowed to take both positive and negative
values with $\sqrt{X_t(x)}$  replaced by $\sigma(t,x,X_t(x))$
($\sigma(\cdot,\cdot,u)$ is H\"older continuous in $u$ of index
$\beta_0>3/4$) in \eqref{1.01}, the pathwise uniqueness
was proved by Mytnik and Perkins \cite{MyP11} and further
investigated by Mytnik and Neuman \cite{MN15}.
Recently, some negative results were obtained. When $\sqrt{X_t(x)}$
is replaced by $|X_t(x)|^{\beta_1}$ in \eqref{1.01}, nonuniqueness
results were obtained for $0<\beta_1<3/4$ in \cite{BMP10,MMyP2014}.
It was also shown in Chen  \cite{ChenY15} that the solution is a
super-Brownian motion with immigration and the pathwise
nonuniqueness holds when a positive function is added on the
right-hand side of \eqref{1.01}. We refer to Li \cite{Li11} and
Xiong \cite{Xiong13book} for introductions on superprocesses  and the related SPDEs.

A novel  approach for studying the strong  uniqueness of
\eqref{1.01} was proposed by Xiong \cite{X13}, where an SPDE for the
\textit{distribution-function process} of measure-valued
super-Brownian motion was formulated and the strong existence and
uniqueness for this SPDE were established. Similar stochastic
equations of distribution-function processes for superprocesses were
also established in Dawson and Li \cite{DLi12}. He \textit{et al.}
\cite{HLY14} showed that the distribution-function process of a
one-dimensional super-L\'evy process with general branching
mechanism is the unique strong solution to
another SPDE. The uniqueness of solution to the martingale problem
of a superprocess with interactive immigration mechanism was
first established in Mytnik and Xiong \cite{MX14}, and
then in Xiong and Yang \cite{XiongY} (for a more general process) by
studying the corresponding SPDE for its distribution-function
process.

Let $\mathcal{X}_0(\mbb{R})$ be the Hilbert space consisting of all
measurable functions $f$ on $\mbb{R}$ satisfying
$\int_{\mbb{R}}|f(x)|e^{-|x|}dx<\infty$. Let $D(\mbb{R})$ be the set
of bounded right-continuous nondecreasing functions $f$ on $\mbb{R}$
satisfying $f(-\infty)= 0$. Let $D_c(\mbb{R})$ be the subset of
$D(\mbb{R})$ consisting of continuous functions.

Given a Polish space $E$ and a $\sigma$-finite Borel measure $\pi$
on $E$,  consider the following SPDE:
 \beqlb\label{1.4}
Y_t(y)=Y_0(y)+\frac12\int_0^t\Delta
Y_s(y)ds+\int_0^t\int_{E}G(u,Y_s(y))W(ds,du),\quad Y_0\in
\mathcal{X}_0(\mbb{R}),
 \eeqlb
where $\{W(ds,du):s\ge0,u\in E\}$ is a Gaussian white noise with
intensity $ds\pi(du)$ and $G$ is a Borel function on $E\times
\mbb{R}_+$ $(\mbb{R}_+:=[0,\infty))$. Then an
$\mathcal{X}_0(\mbb{R})$-valued process $\{Y_t(y):t\geq 0,
y\in\IR\}$ is a solution to (\ref{1.4}) if for any $f\in
C^2_0(\mathbb{R})$ (to be define at the end of this
section) with
 $\int_{\mbb{R}}|f(x)|dx<\infty$,
 \beqlb\label{1.5}
\<Y_t, f\>=\<Y_0,
f\>+\frac{1}{2}\int_0^t\<Y_s,f''\>ds+\int_0^t\int_{E}\Big[\int_\mathbb{R}
G(u,Y_s(y))f(y)dy\Big] W(ds,du),~~\mbf{P}\mbox{-a.s.}
 \eeqlb
Xiong \cite{X13} proved that \eqref{1.4} has a unique
strong $\mathcal{X}_0(\mbb{R})$-valued solution if $G$ satisfies
the following conditions: there is a constant $C>0$ so that
 \beqlb\label{1.6}
\int_{E}|G(u,x)|^2\pi(du)\le C(1+x^2),~~x\in\mbb{R}
 \eeqlb
and
 \beqlb\label{1.7}
\int_{E}|G(u,x_1)-G(u,x_2)|^2\pi(du)\le C|x_1-x_2|,~~
x_1,x_2\in\mbb{R}.
 \eeqlb
For $G(u,v)=1_{\{0\le u\le v\}}$ and $G(u,v)=1_{\{0\le u\le v\le1\}}-v$,
the solutions to \eqref{1.4} are the distribution-function processes
of super-Brownian motion and Fleming-Viot process,
respectively.

 In this paper we further improve the results in \cite{X13}. In particular,
we establish a comparison theorem which shall be of independent
interest. Using it we prove that \eqref{1.4} indeed has a
unique strong $D_c(\mbb{R})$-valued solution when $G$
satisfies certain conditions. We then apply the results to show the
well-posedness of martingale problems for an interacting
super-Brownian motion and an interacting Fleming-Voit process. We
summarize the main results in the following.

\textit{Our first main result} of this paper, Theorem \ref{p2.3},
shows that  SPDE \eqref{1.4} has a unique strong
$D_c(\mbb{R})$-valued solution when  $G$ satisfies the following
conditions:
 \beqlb\label{2.17}
\int_E|G(u,0)|\pi(du)=0,
 \eeqlb
 for each $k\ge1$ there is a constant $C_k>0$ so that
 \beqlb\label{2.1}
\int_E|G(u,x)|^2\pi(du)\le C_k,\qquad x\in[0,k]
 \eeqlb
and
 \beqlb\label{2.2}
\int_E|G(u,x_1)-G(u,x_2)|^2\pi(du)\le C_k|x_1-x_2|,\qquad
 x_1,x_2\in[0,k].
 \eeqlb

In \textit{ our second main results} of this paper, Theorems
\ref{t3.1} and \ref{t3.3}, we prove that the martingale problems for
both the interacting super-Brownian motion and the interacting
Fleming-Viot process are well-posed by first associating the
martingale problems with the corresponding SPDEs (\ref{1.4}) with
$E=\IR_+, \, G(u,x)=1_{\{u\le x\}}\sqrt{\sigma(u)}$ and with
$E=[0,1]^2, \, G(a,b,x)=1_{\{a\leq x\leq b\}}\sqrt{\gamma(a,b)}$ for
nonnegative functions $\sigma$ and $\gamma$ (see Lemmas \ref{l3.1}
and \ref{t4.2}), respectively, and then applying Theorem \ref{p2.3}.
They partially generalize the recent work of
\cite{MX14,XiongY}.

We want to point out  that, in our paper, the existence of solutions
to the martingale problems follows directly from
the relationship with their corresponding SPDEs and the existence of
solutions to these SPDEs, which is different from the
 approach of approximating  martingale problem for the classical superprocesses in
\cite[Proposition 2.4]{MX14} and the approach of approximating
branching interacting particle system in Xiong and Yang
\cite[Theorem 3.3]{XiongY}.

\textit{Our last main results}, Theorems \ref{t3.2}--\ref{t3.6} and
\ref{t4.1}, concern  properties such as the existence of density
fields and survival-extinction behaviors of the
interacting super-Brownian motions and Fleming-Viot processes. In
particular, using martingale arguments we discuss the
survival-extinction properties of the interacting
super-Brownian motion $\{X_t:t\ge0\}$ determined by the local
martingale problem (\ref{MP1a})--(\ref{MP1b}). We show that the total
mass process $\{X_t(1):t\ge0\}$ satisfies $X_t(1)=0$ for all $t$
large enough if $\int_0^x \sigma(y)dy$  decreases to $0$ slow enough
as $x\goto 0+$ and $X_t(1)>0$ for all $t$ if $\int_0^x \sigma(y)dy$
decreases to $0$ fast enough.

The rest of the paper is organized as follows. In Section 2, the
comparison theorem and the strong uniqueness of
$D_c(\mbb{R})$-valued solution to \eqref{1.4} are established. The
interacting super-Brownian motions and Fleming-Viot processes are
studied in Sections 3 and 4, respectively.

{\bf Notations:} We always assume that all random elements are
defined on a filtered complete probability space $(\Omega,
\mathscr{F}, ( \mathscr{F}_t)_{t\ge0},\mathbf{P})$ satisfying the
usual hypotheses. For a topological space $V$ let $\mathscr{B}(V)$
be the space of Borel measurable functions and Borel sets on $V$.
Let $B(V)$ be the collection of bounded functions on $V$ furnished
with the supremum norm $\|\cdot\|$ and $C(V)$ the space of bounded
continuous functions on $V$. For $f,g\in B(\mbb{R})$ write $\<f,g\>=
\int_{\mbb{R}}f(x)g(x)dx$ whenever it exists. Let $C^2(\mathbb{R})$
be the space of twice continuously differentiable functions on
$\mathbb{R}$. Let $C_0(\mathbb{R})$ be the subset of $C(\IR)$
consisting of functions vanishing at infinity and let
$C_0^k(\mathbb{R})$ ($k\ge1$) be the subset of functions with
derivatives up to order $k$ belonging to $C_0(\mathbb{R})$. We use
the superscript ``+'' to denote the subsets of non-negative elements
of the function spaces. Denote $M(\mbb{R})$ for the space of finite
Borel measures on $\mbb{R}$ equipped with the topology of weak
convergence. Then there is an obvious one-to-one correspondence
between
 $M(\mbb{R})$ and $D(\mbb{R})$ by associating a measure with its
 distribution function. We endow $D(\mbb{R})$ with the topology induced by this
correspondence from the weak convergence topology on $M(\mbb{R})$.
Then for any $M(\mbb{R})$-valued stochastic process $\{X_t: t\ge
0\}$, its \textit{distribution-function process} $\{Y_t: t\ge 0\}$
is a $D(\mbb{R})$-valued stochastic process. For $\mu\in M(\mbb{R})$
and $f\in B(\mbb{R})$ write $\mu(f) \equiv \<\mu,f\>:= \int_\mbb{R}
f(x)\mu(dx)$. Write $f(\infty):=\lim_{x\to\infty}f(x)$ and
$f(-\infty):=\lim_{x\to-\infty}f(x)$ for $f\in D(\mbb{R})$. Let
$\lambda$ denote the Lebesgue measure on $\mbb{R}$. Throughout this
paper we use $C$ to denote a positive constant whose value might
change from line to line. We also write $C_n$ if the constant
depends on $n\ge1$.

\section{Existence and strong uniqueness of solution to SPDE}

\setcounter{equation}{0}

In this section we establish the existence and strong uniqueness of
$D(\mbb{R})$-valued solutions to the SPDE \eqref{1.4} under
conditions (\ref{1.6})--(\ref{2.17}) and conditions
(\ref{2.17})--(\ref{2.2}), respectively.
 Xiong \cite{X13} proved that \eqref{1.4}
has a unique strong $\mathcal{X}_0(\mbb{R})$-valued solution under
conditions (\ref{1.6})--(\ref{1.7}). So we only need  to show that
the $\mathcal{X}_0(\mbb{R})$-valued solution is in fact
$D(\mbb{R})$-valued under additional condition (\ref{2.17}) and
further prove that (\ref{1.4}) has a unique strong
$D_c(\mbb{R})$-valued solution under conditions
(\ref{2.17})--(\ref{2.2}).

For this purpose, we  first establish the following
\textit{comparison theorem} for (\ref{1.4}) under the two sets of
conditions, which is of independent interest. The corresponding
strong uniqueness then follows.

 \bproposition\label{l2.1} (i) Let $\{\bar{Y}_{0,t}:t\ge0\}$ and
$\{\bar{\bar{Y}}_{0,t}:t\ge0\}$ be any two $\mathcal{X}_0(\mbb{R})$-valued
solutions of (\ref{1.4}) under conditions (\ref{1.6})--(\ref{1.7}).
If $\bar{Y}_{0,0}(x)\le\bar{\bar{Y}}_{0,0}(x)$ $\mbf{P}$-a.s. for
$\lambda$-a.e. $x$  and
 \beqlb\label{2.3}
\sup_{t\in[0,T]}\Big[\int_{\mbb{R}}\bar{Y}_{0,t}(x)^2e^{-|x|}dx
+\int_{\mbb{R}}\bar{\bar{Y}}_{0,t}(x)^2e^{-|x|}dx\Big]<\infty
\qquad \mbf{P}\mbox{-a.s.}
 \eeqlb
for all $T>0$, then for any $t>0$,
$$\mbf{P}\{\bar{Y}_{0,t}(x)\le
\bar{\bar{Y}}_{0,t}(x) \text{\,\, for \,\,} \lambda \text{-a.e.} \,\, x\}=1.$$

(ii) Let $\{\bar{Y}_{1,t}:t\ge0\}$ and
$\{\bar{\bar{Y}}_{1,t}:t\ge0\}$ be any two $D(\mbb{R})$-valued
solutions of (\ref{1.4}) under conditions (\ref{2.1})--(\ref{2.2}).
If $\bar{Y}_{1,0}(x)\le\bar{\bar{Y}}_{1,0}(x)$ $\mbf{P}$-a.s. for
all $x\in\mbb{R}$ and
 \beqlb\label{2.4}
\sup_{t\in[0,T]}[\bar{Y}_{1,t}(\infty)+\bar{\bar{Y}}_{1,t}(\infty)]<\infty
\qquad \mbf{P}\mbox{-a.s.}
 \eeqlb
for all $T>0$, then for any $t>0$,
$$\mbf{P}\{\bar{Y}_{1,t}(x)\le
\bar{\bar{Y}}_{1,t}(x) \,\,\text{for all\,\, } x\in\mbb{R}\}=1.$$
 \eproposition
 \proof
The proof is inspired by that of \cite[Lemma 4.10]{XiongY} and
\cite[Proposition 3.1]{MX14}. It is divided into four steps.

{\bf Step 1.} For $n\ge1$ and $x,y\in\mbb{R}$ let
 \beqnn
g_n^x(y):=g_n(x-y):=\sqrt{\frac{n}{2\pi}}\exp\Big\{-\frac{n(x-y)^2}{2}\Big\},\quad
\phi_n(x):=\int_{-\infty}^\infty[(x-y)\vee0] g_n(y)dy.
 \eeqnn
Then
$$\phi_n'(x)=\int_{-\infty}^xg_n(y)dy, \,
\phi_n''(x)=g_n(x)\ge0, \, \phi_n(x)\to x\vee0$$ as $n\to\infty$ and
$|x|\phi_n''(x)\le1$ for all $n\ge1$ and $x\in\mbb{R}$. Define
$$J(x):=\int_{\mbb{R}}e^{-|y|}\rho_0(x-y)dy$$ with the mollifier
$\rho_0$ given by
 \beqnn
\rho_0(x):=\tilde{C}\exp\big\{-1/(1-x^2)\big\}1_{\{|x|<1\}},
 \eeqnn
where $\tilde{C}$ is a constant so that
$\int_{\mbb{R}}\rho_0(x)dx=1$. By (2.1) in \cite{M85}, for each
$n\ge0$ there exist constants $\tilde{c}_n,\tilde{C}_n>0$ so that
the $n$th-derivative $J^{(n)}$ of $J$ satisfies
 \beqlb\label{2.5}
\tilde{c}_ne^{-|x|}\le |J^{(n)}(x)|\le \tilde{C}_ne^{-|x|}, \qquad
x\in\mbb{R},
 \eeqlb
which implies
 \beqlb\label{2.6}
|J^{(n)}(x)|\le \tilde{C}_nJ(x), \qquad x\in\mbb{R}.
 \eeqlb

{\bf Step 2.} Let $t>0$ be fixed in the following. For $k\ge1$
define stopping times
 \beqnn
\tau_{0,k}:=\inf\Big\{s\ge0:\<\bar{Y}_{0,s}^2,J\>
+\<\bar{\bar{Y}}_{0,s}^2,J\>>k\Big\}
 \eeqnn
and
 \beqnn
\tau_{1,k}:=\inf\big\{s\ge0:\bar{Y}_{1,s}(\infty)+\bar{\bar{Y}}_{1,s}(\infty)>k\big\}
 \eeqnn
with the convention $\inf\emptyset=\infty$. It follows from
\eqref{2.3} and \eqref{2.4} that
 \beqlb\label{2.7}
\lim_{k\to\infty}\tau_{0,k}=\infty \quad\mbox{and}\quad
\lim_{k\to\infty}\tau_{1,k}=\infty \qquad \mbf{P}\mbox{-a.s}.
 \eeqlb
For $x\in\mbb{R}, m\ge1$ and $i=0,1$, $s\ge0$ let
 \beqnn
Y_{i,s}:=\bar{Y}_{i,s}-\bar{\bar{Y}}_{i,s},
~\tilde{G}_i(s,u,x):=G(u,\bar{Y}_{i,s}(x))-G(u,\bar{\bar{Y}}_{i,s}(x)),
~v_{i,s}^m(x):=\<Y_{i,s},g_m^x\>.
 \eeqnn
By conditioning we may assume that $\bar{Y}_{i,0}$ and
$\bar{\bar{Y}}_{i,0}$ are deterministic. It follows from \eqref{1.5}
and It\^o's formula that for each $x\in\mbb{R}$,
 \beqnn
\phi_n(v_{i,t\wedge\tau_{i,k}}^m(x))
 \ar=\ar
\phi_n(v_{i,0}^m(x))
+\frac12\int_0^{t\wedge\tau_{i,k}}\phi_n'(v_{i,s}^m(x)) \Delta
v_{i,s}^m(x)ds \cr
 \ar\ar
+\frac12\int_0^{t\wedge\tau_{i,k}}\phi_n''(v_{i,s}^m(x))ds
\int_E\Big[\int_{\mbb{R}}\tilde{G}_i(s,u,y)g_m^x(y)dy\Big]^2\pi(du)
+\mbox{mart}.
 \eeqnn
Therefore,
 \beqlb\label{2.8}
 \ar\ar
\mbf{E}\Big\{\int_{\mbb{R}}\phi_n(v_{i,t\wedge\tau_{i,k}}^m(x))J(x)dx\Big\}
\cr
 \ar\ar\quad=
\int_{\mbb{R}}\phi_n(v_{i,0}^m(x))J(x)dx
+\frac12\mbf{E}\Big\{\int_0^{t\wedge\tau_{i,k}}ds\int_{\mbb{R}}\phi_n'(v_{i,s}^m(x))
\Delta v_{i,s}^m(x)J(x)dx\Big\} \cr
 \ar\ar\qquad
+\frac12\mbf{E}\Big\{\int_0^{t\wedge\tau_{i,k}}ds\int_{\mbb{R}}\phi_n''(v_{i,s}^m(x))
J(x)dx\int_E\Big[\int_{\mbb{R}}\tilde{G}_i(s,u,y)g_m^x(y)dy\Big]^2\pi(du)\Big\}
\cr
 \ar\ar\quad=:
I_0(i,m,n)+I_1(i,m,n,k,t)+I_2(i,m,n,k,t).
 \eeqlb

{\bf Step 3.} Observe that
 \beqnn
\int_{\mbb{R}}\phi_n'(v_{i,s}^m(x)) \Delta v_{i,s}^m(x)J(x)dx
=\int_{\mbb{R}}(\phi_n(v_{i,s}^m(x)))'' J(x)dx
-\int_{\mbb{R}}\phi_n''(v_{i,s}^m(x)) |\nabla v_{i,s}^m(x)|^2J(x)dx.
 \eeqnn
Since $\phi_n''=g_n\ge0$, the second term on the right-hand side is
negative. Thus
 \beqlb\label{2.9}
2I_1(i,m,n,k,t)
 \ar\le\ar
\mbf{E}\Big\{\int_0^{t\wedge\tau_{i,k}}ds
\int_{\mbb{R}}(\phi_n(v_{i,s}^m(x)))'' J(x)dx\Big\} \cr
 \ar=\ar
\mbf{E}\Big\{\int_0^{t\wedge\tau_{i,k}}ds
\int_{\mbb{R}}\phi_n(v_{i,s}^m(x))J''(x)dx\Big\} \cr
 \ar\le\ar
C\,\mbf{E}\Big\{\int_0^{t\wedge\tau_{i,k}}ds
\int_{\mbb{R}}\phi_n(v_{i,s}^m(x))J(x)dx\Big\},
 \eeqlb
where integration by parts and \eqref{2.6} are used. By the H\"older
inequality and conditions \eqref{1.7} and \eqref{2.2},
 \beqlb\label{2.10}
2I_2(i,m,n,k,t) \ar\le\ar
\mbf{E}\Big\{\int_0^{t\wedge\tau_{i,k}}ds\int_{\mbb{R}}\phi_n''(v_{i,s}^m(x))
J(x)dx\int_{\mbb{R}}g_m^x(y)dy\int_E\tilde{G}_i(s,u,y)^2\pi(du)\Big\}
\cr
 \ar\le\ar
C_k\mbf{E}\Big\{\int_0^{t\wedge\tau_{i,k}}ds\int_{\mbb{R}}\phi_n''(v_{i,s}^m(x))
\<|Y_{i,s}|,g_m^x\>J(x)dx\Big\}.
 \eeqlb

 For $i=0$, by the fact $\|\phi_n'\|\le1$, $\|\phi_n''\|+\|\phi_n'''\|\le
C_n$ and \cite[Lemma 2.1]{X13}, it is elementary to see that for
$0\leq s\le\tau_{0,k}$,
 \beqlb\label{2.11} \<v_{0,s}^m,J\>\to
\<Y_{0,s},J\>,~ \<|\phi_n(v_{0,s}^m)-\phi_n(Y_{0,s})|,J\> \le
\<|v_{0,s}^m-Y_{0,s}|,J\> \to0,
 \eeqlb
and
 \beqlb\label{2.12}
 \ar\ar
\Big|\int_{\mbb{R}}\phi_n''(v_{0,s}^m(x))\<|Y_{0,s}|,g_m^x\>J(x)dx
-\int_{\mbb{R}}\phi_n''(Y_{0,s}(x))|Y_{0,s}(x)|J(x)dx\Big| \cr
 \ar\le\ar
\int_{\mbb{R}}\Big[\|\phi_n''\|\big|\<|Y_{0,s}|,g_m^x\>-|Y_{0,s}(x)|\big|
+\|\phi_n'''\|\big|v_{0,s}^m(x)-Y_{0,s}(x)\big||Y_{0,s}(x)|\Big]J(x)dx
\cr
 \ar\le\ar
\|\phi_n''\|\int_{\mbb{R}}\Big[\big|\<|Y_{0,s}|,g_m^x\>-|Y_{0,s}(x)|\big|J(x)dx
\cr
 \ar\ar
+\|\phi_n'''\|\Big[\int_{\mbb{R}}\big|v_{0,s}^m(x)-Y_{0,s}(x)\big|^2J(x)dx
\cdot\int_{\mbb{R}}|Y_{0,s}(x)|^2J(x)dx\Big]^{\frac12} \to0
 \eeqlb
as $m\to\infty$. Observe that
 \beqlb\label{2.13}
\phi_n(x)\le |x|+\int_{\mbb{R}}|y|g_n(y)dy \le|x|+1
 \eeqlb
for all $n\ge1$. Then by (2.13) in \cite{X13},
 \beqlb\label{2.14}
\sup_{m\ge1}\Big|\int_{\mbb{R}}\phi_n(v_{0,s}^m(x))J(x)dx\Big|
+\Big|\int_{\mbb{R}}\phi_n''(v_{0,s}^m(x))\<|Y_{0,s}|,g_m^x\>J(x)dx\Big|
< \infty,~\mbox{on~}\{s\le \tau_{0,k}\}.
 \eeqlb
By \cite[Lemma 4.5]{XiongY}, one can check that
\eqref{2.11}--\eqref{2.12} and \eqref{2.14} also hold for $i=1$.
Combining this with \eqref{2.9}--\eqref{2.10} and dominated
convergence we get
 \beqnn
\lim_{m\to\infty}\mbf{E}\Big\{\int_{\mbb{R}}
\phi_n(v_{i,t\wedge\tau_{i,k}}^m(x))J(x)dx\Big\}
=\mbf{E}\Big\{\int_{\mbb{R}}\phi_n(Y_{i,t\wedge\tau_{i,k}}(x))J(x)dx\Big\}
 \eeqnn
and
 \beqnn
 \ar\ar
\lim_{m\to\infty}I_0(i,m,n)=\<\phi_n(Y_{i,0}),J\>=:J_0(i,n),\cr
 \ar\ar
\limsup_{m\to\infty}I_1(i,m,n,k,t)\le
C\,\mbf{E}\Big\{\int_0^{t\wedge\tau_{i,k}}ds
\int_{\mbb{R}}\phi_n(Y_{i,s}(x))J(x)dx\Big\} =:J_1(i,n,k,t), \cr
 \ar\ar
\limsup_{m\to\infty}I_2(i,m,n,k,t) \le
C_k\mbf{E}\Big\{\int_0^{t\wedge\tau_{i,k}}ds\int_{\mbb{R}}\phi_n''(Y_{i,s}(x))
|Y_{i,s}(x)|J(x)dx\Big\}=:J_2(i,n,k,t)
 \eeqnn
for both $i=0,1$.

Now taking $m\to\infty$ in \eqref{2.8} we have
 \beqlb\label{2.15}
\mbf{E}\Big\{\int_{\mbb{R}}\phi_n(Y_{i,t\wedge\tau_{i,k}}(x))J(x)dx\Big\}
\le J_0(i,n)+J_1(i,n,k,t)+J_2(i,n,k,t).
 \eeqlb

{\bf Step 4.} Since $0\le|x|\phi_n''(x)\le1$ for all
$x\in\mbb{R},n\ge1$ and $\lim_{n\to\infty}|x|\phi_n''(x)=0$ for each
$x\in\mbb{R}$, by dominated convergence we have
 \beqlb\label{2.16}
\lim_{n\to\infty}J_2(i,n,k,t)\le
C_k\mbf{E}\Big\{\int_0^{t\wedge\tau_{i,k}}ds
\int_{\mbb{R}}\lim_{n\to\infty}\phi_n''(Y_{i,s}(x))
|Y_{i,s}(x)|J(x)dx\Big\}
 =0.
 \eeqlb
Using dominated convergence and \eqref{2.13} we know that
 \beqnn
 \ar\ar
\lim_{n\to\infty}J_0(i,n)=\int_{\mbb{R}}(Y_{i,0}(x))^+J(x)dx=0, \cr
 \ar\ar
\lim_{n\to\infty}J_1(i,n,k)=C\,\mbf{E}\Big\{\int_0^{t\wedge\tau_{i,k}}ds
\int_{\mbb{R}}(Y_{i,s}(x))^+J(x)dx\Big\},\cr
 \ar\ar
\lim_{n\to\infty}\mbf{E}\Big\{\int_{\mbb{R}}\phi_n(Y_{i,t\wedge\tau_{i,k}}(x))J(x)dx\Big\}
=\mbf{E}\Big\{\int_{\mbb{R}}\Big(Y_{i,t\wedge\tau_{i,k}}(x)\Big)^+J(x)dx\Big\}.
 \eeqnn
Combining with \eqref{2.15}--\eqref{2.16} we get
 \beqnn
\mbf{E}\Big\{1_{\{t\le\tau_{i,k}\}}\int_{\mbb{R}}(Y_{i,t}(x))^+J(x)dx\Big\}
\ar\le\ar
\mbf{E}\Big\{\int_{\mbb{R}}\Big(Y_{i,t\wedge\tau_{i,k}}(x)\Big)^+J(x)dx\Big\}\cr
 \ar\le\ar
 C\,\mbf{E}\Big\{\int_0^{t\wedge\tau_{i,k}}ds
\int_{\mbb{R}}(Y_{i,s}(x))^+J(x)dx\Big\}\cr \ar=\ar
C\,\mbf{E}\Big\{\int_0^t\Big[1_{\{s\le\tau_{i,k}\}}
\int_{\mbb{R}}(Y_{i,s}(x))^+J(x)dx\Big]ds\Big\}.
 \eeqnn
Then by Gronwall's lemma,
 \beqnn
\mbf{E}\Big\{1_{\{t\le\tau_{i,k}\}}\int_{\mbb{R}}(Y_{i,t}(x))^+J(x)dx\Big\}=0,
 \eeqnn
which implies that
 \beqnn
1_{\{t\le\tau_{i,k}\}}\int_{\mbb{R}}(Y_{i,t}(x))^+J(x)dx=0,\quad
\mbf{P}\mbox{-a.s}.
 \eeqnn
Letting $k\to\infty$ we see that
$\int_{\mbb{R}}(Y_{i,t}(x))^+J(x)dx=0$ $\mbf{P}$-a.s. by
\eqref{2.7}. Then the desired result follows.
 \qed

 \blemma\label{l2.2}
 Suppose that conditions \eqref{1.6}--\eqref{2.17}
hold and  $Y=\{Y_t:t\ge0\}$ is a solution of \eqref{1.4} satisfying
$\sup_{t\in[0,T]}\mbf{E}\{\< Y_t^2,e^{-|\cdot|}\>\}<\infty$ for all
$T>0$. If $Y_0\in D(\mbb{R})$, then $Y$ has a $
C([0,\infty)\times\mbb{R})$-valued modification
$\tilde{Y}=\{\tilde{Y}_t:t\ge0\}$ with $\tilde{Y}_0\in D(\IR)$ (i.e.
$\mbf{P}\{\<|Y_t-\tilde{Y}_t|,1\>=0\}=1$ for all $t\ge0$). Moreover,
$\mbf{P}\{\tilde{Y}_t\in D_c(\mbb{R})$ for all $t>0\}=1$ and
$\sup_{t\in[0,T]}\tilde{Y}_t(\infty)<\infty$, $\mbf{P}$-a.s.
 \elemma
 \proof
The proof is proceeded in four steps.

{\bf Step 1.} By the assertion and proof of \cite[Theorem 1.2]{X13},
the solution $Y$ of \eqref{1.4} has a $
C([0,\infty)\times\mbb{R})$-valued modification
$\tilde{Y}=\{\tilde{Y}_t:t\ge0\}$. Observe that $Y_t\equiv0$ is also
a solution of \eqref{1.4}. It then follows from Proposition
\ref{l2.1}(i) that $\mbf{P}\{\tilde{Y}_t(x)\ge0$ for all $t\ge0$ and
$x\in\mbb{R}\}=1$.

Fixed $y_1\ge y_2$, let $Y_{1,t}(x):=Y_t(x+y_1)$  and
$Y_{2,t}(x):=Y_t(x+y_2)$. It follows from Proposition \ref{l2.1}(i)
again that $Y_{1,0}(x)\ge Y_{2,0}(x)$ for all $x\in\mbb{R}$ and
 \beqnn
\mbf{P}\{Y_t(x+y_1)\ge Y_t(x+y_2)~\text{for}~
\lambda\text{-a.e.}\,x\}=1
 \eeqnn
which implies
 \beqnn
\mbf{P}\{\tilde{Y}_t(x)\ge \tilde{Y}_t(y)~\text{for all}~x\ge
y\in\IR\}=1.
 \eeqnn

{\bf Step 2.} It follows from (2.12) in \cite{X13} that, for each
$|f|\le CJ$, we have for any $t>0$,
 \beqnn
\<Y_t,f\>=\<Y_0,P_tf\>+\int_0^t\int_E
\Big[\int_{\mbb{R}}G(u,Y_s(z))P_{t-s}f(z)dz\Big]W(ds,du)\quad
\mbf{P}\mbox{-a.s}.
 \eeqnn
Without loss of generality we assume that $Y_0$ is deterministic in
the following. Thus $\mbf{E}\{\<Y_t,J_n\>\}=\<Y_0,P_tJ_n\>$ for
$J_n(x):=J(x-n)$. By monotone convergence and $\<1,J\>=2$ we have
 \beqnn
\lim_{n\to\infty}\mbf{E}\{\<Y_t,J_n\>\}=
\lim_{n\to\infty}\mbf{E}\Big\{\int_{\mbb{R}}\tilde{Y}_t(x+n)J(x)\Big\}
=\<1,J\>\mbf{E}\{\tilde{Y}_t(\infty)\}=2\mbf{E}\{\tilde{Y}_t(\infty)\}
 \eeqnn
and
 \beqnn
\lim_{n\to\infty}\<Y_0,P_tJ_n\>
 =
\lim_{n\to\infty}\int_{\mbb{R}}Y_0(z+n)dz\int_{\mbb{R}}p_t(y)J(z-y)dy
=2Y_0(\infty),
 \eeqnn
which implies
 \beqlb\label{2.18}
\mbf{E}\{\tilde{Y}_t(\infty)\}=Y_0(\infty).
 \eeqlb
Similarly, by dominated convergence and $\<1,J_n\>=2$ we also get
 \beqlb\label{2.19}
2\mbf{E}\{\tilde{Y}_t(-\infty)\}=\mbf{E}[\lim_{n\to-\infty}\<Y_t,J_n\>]
=\lim_{n\to-\infty}\mbf{E}[\<Y_t,J_n\>]
=\lim_{n\to-\infty}\<Y_0,P_tJ_n\>=0.
 \eeqlb

{\bf Step 3.} In this step we show that
$\sup_{t\in[0,T]}\tilde{Y}_t(\infty)<\infty$ for each $T>0$. By
\eqref{1.5} and the Burkholder-Davis-Gundy inequality,
 \beqlb\label{2.20}
\mbf{E}\Big\{\sup_{t\in[0,T]}\<Y_t,J_n\>\Big\}
 \ar\le\ar
\<Y_0,J_n\> +\frac12\int_0^T\mbf{E}\{\<Y_s,|J_n''|\>\}ds \cr
 \ar\ar
+C\mbf{E}\Big\{\Big|\int_0^Tds\int_E\<G(u,Y_s),J_n\>^2\pi(du)\Big|^{\frac12}\Big\}.
 \eeqlb
Under conditions \eqref{1.7}--\eqref{2.17}, for each $y\ge0$,
 \beqnn
\int_EG(u,y)^2\pi(du) =\int_E|G(u,y)-G(u,0)|^2\pi(du) \le C|y|.
 \eeqnn
It then follows from the H\"older inequality and $\<1,J_n\>=2$ that
 \beqnn
 \ar\ar
\int_E\<G(u,Y_s),J_n\>^2\pi(du)
 \le
2\int_E\<G(u,Y_s)^2,J_n\>\pi(du) \cr
 \ar=\ar
2\int_{\mbb{R}}J_n(x)dx\int_EG(u,Y_s(x))^2\pi(du)
 \le
C\int_{\mbb{R}}|Y_s(x)|J_n(x)dx.
 \eeqnn
Putting together \eqref{2.20} and \eqref{2.6} we have
 \beqlb\label{2.21}
\mbf{E}\Big\{\sup_{t\in[0,T]}\<Y_t,J_n\>\Big\}
 \ar\le\ar
\int_{\mbb{R}}Y_0(x+n)J(x)dx +
C\int_0^Tdt\int_{\mbb{R}}\mbf{E}\{Y_t(x+n)\}J(x)dx \cr
 \ar\ar
+C\Big|\int_0^Tdt\int_{\mbb{R}}\mbf{E}\{Y_t(x+n)\}J(x)dx\Big|^{\frac12}.
 \eeqlb
It then follows from \eqref{2.18} and \eqref{2.21} that for each
$T>0$,
 \beqnn
\sup_{n\ge1}\mbf{E}\Big\{\sup_{t\in[0,T]}\<Y_t,J_n\>\Big\} \le
2Y_0(\infty)+CTY_0(\infty)+C[TY_0(\infty)]^{\frac12}.
 \eeqnn
Thus, by monotone convergence
 \beqlb\label{2.22}
 \ar\ar
\<1,J\>\mbf{E}\Big\{\sup_{t\in[0,T]}\tilde{Y}_t(\infty)\Big\}
 =
\mbf{E}\Big\{\sup_{t\in[0,T],n\ge1}\int_{\mbb{R}}\tilde{Y}_t(x+n)J(x)dx\Big\}
\cr
 \ar=\ar
\mbf{E}\Big\{\lim_{n\to\infty}\sup_{t\in[0,T]}\int_{\mbb{R}}\tilde{Y}_t(x+n)J(x)dx\Big\}
 =
\lim_{n\to\infty}\mbf{E}\Big\{\sup_{t\in[0,T]}\<Y_t,J_n\>\Big\},
 \eeqlb
which implies that $\sup_{t\in[0,T]}\tilde{Y}_t(\infty)<\infty$,
$\mbf{P}$-a.s.

{\bf Step 4.} By \eqref{2.19}, \eqref{2.21} and dominated
convergence
 \beqnn
\ar\ar
\<1,J\>\mbf{E}\Big\{\sup_{t\in[0,T]}\tilde{Y}_t(-\infty)\Big\}
 \le
\mbf{E}\Big\{\sup_{t\in[0,T]}\int_{\mbb{R}}\tilde{Y}_t(-\infty)J(x)dx\Big\}\cr
 \ar\le\ar
\mbf{E}\Big\{\lim_{n\to-\infty}\sup_{t\in[0,T]}\int_{\mbb{R}}\tilde{Y}_t(x+n)J(x)dx\Big\}
 \le
\lim_{n\to-\infty}\mbf{E}\Big\{\sup_{t\in[0,T]}\<\tilde{Y}_t,J_n\>\Big\}
=0,
 \eeqnn
which implies $\sup_{t\in[0,T]}\tilde{Y}_t(-\infty)=0$,
$\mbf{P}$-a.s. Therefore, $\mbf{P}\{\tilde{Y}_t\in D(\mbb{R})\}=1$
for each $t>0$, which finishes the proof. \qed

Now we are ready to present the main theorem.

 \btheorem\label{p2.3}
Suppose that conditions \eqref{2.17}--\eqref{2.2} hold. Then for any
$Y_0\in D(\IR)$, SPDE \eqref{1.4} has a pathwise unique continuous
$D_c(\mbb{R})$-valued solution $\{Y_t:t>0\}$.
 \etheorem

\noindent{\it Proof.} The pathwise uniqueness follows immediately
from Proposition \ref{l2.1}(ii). We prove the existence of a
$D_c(\mbb{R})$-valued solution to \eqref{1.4} in the following. We
assume that $Y_0$ is deterministic. By \cite[Theorem 1.2]{X13} and
Lemma \ref{l2.2}, the following SPDE has a unique strong continuous
$D_c(\mbb{R})$-valued solution:
 \beqlb\label{2.23}
Y_t(y)=Y_0(y)+\frac{1}{2}\int_0^t\Delta
Y_s(y)ds+\int_0^t\int_EG_m(u,Y_s(y)) W(ds,du),
 \eeqlb
where $G_m(u,y):=G(u,y\wedge m)$ for $m\ge1$. For each $m\ge1$, let
$\{Y_t^m:t\ge0\}$ be a strong continuous $D_c(\mbb{R})$-valued
solution to (\ref{2.23}) with $Y_0^m:=Y_0$. For each $m,n\ge1$,
define stopping time $\tau_n^m:=\inf\{t\ge0:Y_t^m(\infty)\ge
n\}$. It then follows from the uniqueness of strong solution to (\ref{2.23}) that
$\mbf{P}$-a.s.
 \beqnn
Y_{t\wedge\tau_m^m\wedge\tau_m^{m+1}}^m=Y_{t\wedge\tau_m^m\wedge\tau_m^{m+1}}^{m+1},
\qquad t\ge0.
 \eeqnn
If $\tau^m_m\le \tau^{m+1}_m$, then
$Y^{m+1}_{\tau^m_m}(\infty)=Y^m_{\tau^m_m}(\infty)= m$, which, by
the definition of $\tau^{m+1}_m$, implies that
$\tau^m_m\geq\tau^{m+1}_m $. Therefore, $\tau^m_m=\tau^{m+1}_m $ and
$Y^m_{\cdot\wedge \tau^m_m}=Y^{m+1}_{\cdot\wedge \tau^{m+1}_m} $.
This implies that $\tau_m^m$ is increasing in $m$. Observe that for
$J_n(x)=J(x-n)$ we have
$\<Y_t^m,J_n\>=\int_{\mbb{R}}Y_t^m(y+n)J(y)dy$. It then follows from
\eqref{2.23} that
 \beqnn
\int_{\mbb{R}}Y_t^m(y+n)J(y)dy
 \ar=\ar
\int_{\mbb{R}}Y_0(y+n)J(y)dy +\frac{1}{2}\int_0^tds
\int_{\mbb{R}}Y_s^m(y+n)J''(y)dy \cr
 \ar\ar
+\int_0^t\int_E\Big[\int_{\mbb{R}}G_m(u,Y_s^m(y+n))J(y)dy\Big]
W(ds,du).
 \eeqnn
Letting $n\to\infty$, we have
 \beqnn
Y_t^m(\infty)=Y_0(\infty)+\int_0^t\int_EG_m(u,Y_s^m(\infty))
W(ds,du),
 \eeqnn
which implies $\mbf{E}\{Y^m_{\tau_m^m\wedge
t}(\infty)\}=Y_0(\infty)$. Thus
 \beqnn
m\mbf{P}\{\tau_m^m<t\}=
\mbf{E}\big\{Y^m_{\tau_m^m}(\infty)1_{\{\tau_m^m<t\}}\big\}\leq
\mbf{E}\{Y^m_{\tau_m^m\wedge t}(\infty)\}=Y_0(\infty).
 \eeqnn
Letting $t\rightarrow \infty$, we have
 \beqnn
\mbf{P}\{\tau_m^m<\infty\}\leq Y_0(\infty)/m\to0
 \eeqnn
as $m\to\infty$. So, we must have $\tau_m^m\uparrow\infty$.
Therefore, we can define the solution to (\ref{1.4}) by
$Y_t=Y^m_{t\wedge\tau_m^m}$ for $t\le \tau_m^m$. This completes the
proof. \qed

\section{Interacting super-Brownian motions}

\setcounter{equation}{0}

In the early work of M\'{e}l\'{e}ard and Roelly \cite{MR93} a
measure-valued branching process with mean field interaction was
introduced. Using particle system approximation, it was shown that if
$\tilde{\sigma}\in C(M(\mathbb{R})\times\mathbb{R})^+$, then the
continuous $M(\IR)$-valued solution $\{X_t:t\ge0\}$ to the following
martingale problem exists: for any $f\in C^2(\mbb{R})$,
 \begin{equation}\label{MP0a}
M_t^f:=X_t(f)-X_0(f)-\frac12\int_0^t X_s(f'')ds, \qquad t\ge0
 \end{equation}
is a square-integrable continuous martingale with quadratic
variation process
 \begin{equation}\label{MP0b}
\< M^f\>_t=\int_0^t \< X_s, \tilde{\sigma}(X_s)f^2\> ds.
 \end{equation}
However, owing to the interaction, the fundamental multiplicative
property of the measure-valued branching process
disappears and one cannot associate  the interacting
measure-valued branching processes with a cumulant
semigroup, which can be characterized as the unique
solution of a non-linear partial differential equation. The
uniqueness of solution to the above martingale problem is still
unknown.  Using snake representation and random time change
techniques, Delmas and Dhersin \cite{DD03} proved the existence of
solution to a similar martingale problem for interacting
super-Brownian motion. But again, the uniqueness of solution is left
open in general.

In this section, we consider an \textit{interacting super-Brownian motion}
$\{X_t:t\ge0\}$, which is a continuous $M(\IR)$-valued process solving the
following local martingale problem: for any $f\in C^2(\mbb{R})$,
 \begin{equation}\label{MP1a}
M_t(f):=X_t(f)-X_0(f)-\frac12\int_0^tX_s(f'')ds, \qquad t\ge0
 \end{equation}
is a continuous local martingale with quadratic variation process
 \begin{equation}\label{MP1b}
\<M(f)\>_t=\int_0^tds\int_{\mbb{R}}\sigma(X_s, x)f(x)^2X_s(dx)
:=\int_0^tds\int_{\mbb{R}}\sigma(X_s(-\infty,x]))f(x)^2X_s(dx),
 \end{equation}
where the {\it branching rate function} 
$\sigma\in
\mathscr{B}(\mbb{R_+})^+$ is bounded on $[0,n]$ for each $n\ge1$ and
 \beqnn
\lambda\{u\in\mbb{R}_+:\sigma(u)=0\}=0.
 \eeqnn
Note that the local martingale problem \eqref{MP1a}--\eqref{MP1b} is
not covered by that of \eqref{MP0a}--\eqref{MP0b}. For constant
$\sigma$, super-Brownian motion with branching rate $\sigma$ is the
unique solution to the above-mentioned martingale problem.

We will establish both the existence and uniqueness of solution to martingale
problem (\ref{MP1a})-(\ref{MP1b}) by associating its solution with an SPDE of type (\ref{1.4}). Then properties
of the interacting super-Brownian motion will be investigated.

\subsection{Well-posedness of the martingale problem}

The next result is on the connection between the local martingale
problem and the corresponding distribution-function-valued SPDE. We
say that a $D(\mbb{R})$-valued process $\{Y_t:t\ge0\}$ is a {\it
distribution-function process} of a $M(\IR)$-valued process
$\{X_t:t\ge0\}$ if $Y_t(x)=X_t((-\infty, x]) $ for all $t\geq 0$ and
$x\in\IR$.

\blemma\label{l3.1}
A continuous $D(\mbb{R})$-valued process $\{Y_t:t\ge0\}$ is the
distribution-function process of the interacting
super-Brownian motion if and only if there is, on an enlarged
probability space, a Gaussian white noise $\{W(dt,du):t\ge0,u\ge0\}$
with intensity $dtdu$ so that $\{Y_t:t\ge0\}$ solves the following
SPDE:
 \beqlb\label{3.1}
Y_t(x)=Y_0(x)+\frac12\int_0^t\Delta
Y_s(x)ds+\int_0^t\int_0^\infty1_{\{u\le
Y_s(x)\}}\sqrt{\sigma(u)}W(ds,du).
 \eeqlb
 \elemma

\noindent{\bf Remark }Choose $E=\IR_+$, $\pi=$Lebesgue measure on $\mbb{R}_+$ and
 \beqnn
G(u,x)=1_{\{u\le x\}}\sqrt{\sigma(u)}, \qquad u,x\in\IR_+.
 \eeqnn
Then $\{Y_t:t\ge0\}$, the $D(\IR)$-valued solution to SPDE
(\ref{1.4})  under conditions \eqref{2.17}--\eqref{2.2}, is in fact
the distribution-function process of the measure-valued process
$\{X_t:t\ge0\}$ which solves the local martingale problem
(\ref{MP1a})--(\ref{MP1b}). Consequently, the existence of solution
to the local martingale problem follows from the existence of
solution to (\ref{1.4}), which is different from the approaches in
\cite{DD03,MR93}, and the uniqueness of  solution follows from the
strong uniqueness of the SPDE. The following theorem
partially generalizes the results of \cite{MX14,XiongY} in which
$\sigma$ is required to be bounded.

\btheorem\label{t3.1}
The local martingale problem (\ref{MP1a})--(\ref{MP1b}) is
well-posed.
 \etheorem

\noindent{\it Proof.} The result follows from
Theorem \ref{p2.3} and Lemma \ref{l3.1} immediately.
The uniqueness of solution to the local martingale problem can alternatively be shown by \cite[Theorem 1.2(b)]{MX14} and a standard localization argument.
 \qed

\noindent{\it Proof of Lemma \ref{l3.1}.} The proof is a
modification of those of \cite[Theorem 1.3]{X13} and \cite[Theorem
3.1]{HLY14}.  Suppose that
$\{Y_t:t\ge0\}$ is a continuous $D(\mbb{R})$-valued solution of
\eqref{3.1}. By integration by parts, $X_t(f)=-\<Y_t,f'\>$ for each
$f\in C_0^3(\mbb{R})$. It then follows from \eqref{3.1} that
 \beqnn
X_t(f)
 \ar=\ar
-\<Y_0,f'\>-\frac12\int_0^t\<Y_s,f'''\>ds
-\int_0^t\int_0^\infty\Big[\int_{\mbb{R}}f'(x)1_{\{u\le
Y_s(x)\}}dx\Big]\sqrt{\sigma(u)}W(ds,du) \cr
 \ar=\ar
X_0(f)+\frac12\int_0^tX_s(f'')ds
-\int_0^t\int_0^\infty\Big[\int_{\mbb{R}}f'(x)1_{\{Y_s^{-1}(u)\le
x\}}dx\Big]\sqrt{\sigma(u)}W(ds,du) \cr
 \ar=\ar
X_0(f)+\frac12\int_0^tX_s(f'')ds+\tilde{I}_t(f),
 \eeqnn
where
 \beqlb\label{3.11}
Y_t^{-1}(u):=\inf\{x\in\mbb{R}:Y_t(x)\ge u\}
 \eeqlb
with the convention $\inf\emptyset=\infty $ and
 \beqnn
\tilde{I}_t(f):=\int_0^t\int_0^\infty
f(Y_s^{-1}(u))\sqrt{\sigma(u)}W(ds,du).
 \eeqnn
One can see that $\{\tilde{I}_t(f):t\ge0\}$ is a continuous local
martingale with quadratic variation process
 \beqnn
\<\tilde{I}(f)\>_t = \int_0^t\int_0^\infty
\sigma(u)f(Y_s^{-1}(u))^2du =\int_0^tds\int_{\mbb{R}}
\sigma(X_s(-\infty,x])f(x)^2X_s(dx).
 \eeqnn
By an approximation argument, one can see the above relation remains
true for all $f\in C^2(\mbb{R})$. Therefore, $\{X_t:t\ge0\}$ is an
interacting super-Brownian motion.

Suppose that $\{X_t:t\ge0\}$ is an interacting super-Brownian motion
determined by the local martingale problem
(\ref{MP1a})--(\ref{MP1b}). We will show that the
distribution-function process $\{Y_t:t\ge0\}$ of $\{X_t:t\ge0\}$
solves \eqref{3.1}. We assume that $X_0$ is deterministic in the
following.

Observe that
 \beqnn
X_t(1)=X_0(1)+M_t(1),
 \eeqnn
where $\{M_t(1):t\ge0\}$ is a continuous local martingale and
$\{M_{t\wedge\tau_n}(1):t\ge0\}$ is a continuous martingale for each
$n\ge1$ with stopping time $\tau_n$ defined by
$\tau_n:=\inf\{t\ge0:X_t(1)\ge n\}$.  It then follows that
$\mbf{E}\{X_{t\wedge\tau_n}(1)\}= X_0(1)$. Since $t\mapsto X_t(1)$
is continuous, $\tau_n\to\infty$ as $n\to\infty$. Thus by Fatou's
lemma,
 \beqlb\label{3.2}
\mbf{E}\{X_t(1)\}=\mbf{E}\Big\{\liminf_{n\to\infty}X_{t\wedge\tau_n}(1)\Big\}
\le
\liminf_{n\to\infty}\mbf{E}\big\{X_{t\wedge\tau_n}(1)\big\}=X_0(1).
 \eeqlb
We can complete the proof in four steps.

{\bf Step 1.} We say
$\{\bar{M}_t(B):t\ge0,B\in\mathscr{B}(\mbb{R})\}$ is a
\textit{continuous local martingale measure} if there are stopping
times $\tilde{\sigma}_n$ with $\tilde{\sigma}_n\to\infty$
$\mbf{P}$-a.s. so that for each $n\ge1$,
$\{\bar{M}_{t\wedge\tilde{\sigma}_n}(B):t\ge0,B\in\mathscr{B}(\mbb{R})\}$
is a continuous martingale measure. That is for each
$B\in\mathscr{B}(\mbb{R})$,
$\{\bar{M}_{t\wedge\tilde{\sigma}_n}(B):t\ge0\}$ is a
square-integrable continuous martingale and for every disjoint
sequence $\{B_1,B_2,\cdots\}\subset\mathscr{B}(\mbb{R})$ we have
 \beqnn
\bar{M}_{t\wedge\tilde{\sigma}_n}\Big(\bigcup_{k=1}^\infty
B_k\Big)=\sum_{k=1}^\infty \bar{M}_{t\wedge\tilde{\sigma}_n}(B_k)
 \eeqnn
by the convergence in $L^2(\Omega,\mbf{P})$. For any $n\ge1$ and $f\in
B(\mbb{R})$, let $M_t^n(f)=M_{t\wedge\tau_n}(f)$ which determines a
continuous martingale measure
$\{M_t^n(B):t\ge0,B\in\mathscr{B}(\mbb{R})\}$ by the same argument
as in the proof of \cite[Theorem 7.25]{Li11}. Then by the same
argument as in \cite[Lemma 5.6]{M02} we have
$M_t^n=M_{t\wedge\tau_n}^{n+1}$ for $n\ge1$. Thus we can define the
continuous local martingale measure
$\{M_t(B):t\ge0,B\in\mathscr{B}(\mbb{R})\}$ by $M_t=M_t^n$ for
$t\le\tau_n$.

{\bf Step 2.} Let $\mathbf{B}(\mbb{R}_+\times\mbb{R}\times\Omega)$
be the space of progressively measurable functions on
$\mbb{R}_+\times\mbb{R}\times\Omega$ and
 \beqnn
\mathcal{F}
 \ar=\ar
\Big\{F\in\mathbf{B}(\mbb{R}_+\times\mbb{R}\times\Omega):
\int_0^tds\int_{\mbb{R}}\sigma(X_s(-\infty,x])F(s,x)^2X_s(dx)<\infty,
~t>0,~\mbf{P}\mbox{-a.s.}\Big\}, \cr \mathcal{F}_b
 \ar=\ar
\Big\{F\in\mathbf{B}(\mbb{R}_+\times\mbb{R}\times\Omega):
\mbf{E}\Big\{\int_0^tds\int_{\mbb{R}}\sigma(X_s(-\infty,x])F(s,x)^2X_s(dx)\Big\}<\infty,~\forall
t>0\Big\}.
 \eeqnn
In this step we show that the stochastic integral of $F\in
\mathcal{F}$ with respect to $\{M(ds,dx):s\ge0,x\in\mbb{R}\}$ is
well defined. By the same argument as in \cite[pp.160--166]{Li11},
the stochastic integral of $F\in\mathcal{F}_b$ with respect to
$\{M(ds,dx):s\ge0,x\in\mbb{R}\}$ is well defined. For $F\in
\mathcal{F}$ define stopping times $\tau_n'$ by
 \beqnn
\tau_n'=\inf\Big\{t\ge0:\int_0^tds\int_{\mbb{R}}\sigma(X_s(-\infty,x])F(s,x)^2X_s(dx)\ge
n\Big\}.
 \eeqnn
Then $\tau_n'\to\infty$ $\mbf{P}$-a.s. and the stochastic integral
of $F\in\mathcal{F}$ with respect to
$\{M(ds,dx):s\ge0,x\in\mbb{R}\}$ is defined by
 \beqnn
\int_0^t\int_{\mbb{R}}F(s,x)M(ds,dx)=\int_0^t\int_{\mbb{R}}F(s,x)1_{\{s\le\tau_n'\}}M(ds,dx),
\qquad t\le\tau_n'.
 \eeqnn

{\bf Step 3.}
 For $g\in B(\mbb{R}_+)$ define continuous local martingale $t\mapsto
Z_t(g)$ by
 \beqlb\label{3.3}
Z_t(g)=\int_0^t\int_{\mbb{R}}g(Y_s(x))\sigma(Y_s(x))^{-\frac12}
1_{\{\sigma(Y_s(x))\neq0\}} M(ds,dx).
  \eeqlb
Observe that the quadratic variation process of $t\mapsto Z_t(g)$
satisfies
 \beqnn
\<Z_t(g)\>_t
 \ar=\ar
\int_0^tds\int_{\mbb{R}}g(Y_s(x))^2\sigma(Y_s(x))^{-1}
 1_{\{\sigma(Y_s(x))\neq0\}}\sigma(Y_s(x))X_s(dx) \cr
 \ar=\ar
\int_0^tds\int_{\mbb{R}}g(Y_s(x))^2\sigma(Y_s(x))^{-1}
 1_{\{\sigma(Y_s(x))\neq0\}}dY_s(x)\cr
 \ar=\ar
\int_0^tds\int_0^\infty g(u)^2
 1_{\{\sigma(u)\neq0, u\le Y_s(\infty)\}}du
=
\int_0^tds\int_{\mbb{R}}g(u)^2 1_{\{u\le Y_s(\infty)\}}du,
 \eeqnn
where $\lambda\{u\in\mbb{R}_+:\sigma(u)=0\}=0$ was used in the last
equation. Combining this with \eqref{3.2} one sees that
$\{Z_t(g):t\ge0\}$ is a martingale for each $g\in B(\mbb{R}_+)$.
Then the family $\{Z_t(g) : t\ge0, g\in B(\mbb{R}_+)\}$ determines a
martingale measure $\{Z(dt,dx) : t\ge0, x\ge0\}$. Then by
\cite[Theorem III-6]{{EM90}}, on an extended  probability space,
there is a Gaussian white noise $\{W(ds,du):s\ge0,u>0\}$ so that
 \beqnn
Z_t(g)=\int_0^t\int_0^\infty g(u)1_{\{u\le Y_s(\infty)\}}W(ds,du).
 \eeqnn
Then we can extend the definition of the stochastic integral of
 \beqnn
g\in
\mathcal{F}_0:=\Big\{g\in\mathbf{B}(\mbb{R}_+\times\mbb{R}_+\times\Omega):
\int_0^tds\int_0^\infty g(s,u)^21_{\{u\le Y_s(\infty)\}}du<\infty
~\mbf{P}\mbox{-a.s.},~\forall t\ge0\Big\}
 \eeqnn
with respect to the martingale measure $\{Z(dt,dx): t\ge0, x\ge0\}$
and
 \beqlb\label{3.4}
\int_0^t\int_0^\infty g(s,x)Z(ds,dx) = \int_0^t\int_0^\infty
g(s,u)1_{\{u\le Y_s(\infty)\}}W(ds,du).
 \eeqlb
It follows from Step 2 and \eqref{3.3} that for each $g\in
\mathcal{F}_0$,
 \beqnn
\int_0^t\int_0^\infty
g(s,u)Z(dt,du)=\int_0^t\int_{\mbb{R}}g(s,Y_s(x))\sigma(Y_s(x))^{-\frac12}
1_{\{\sigma(Y_s(x))\neq0\}} M(ds,dx).
 \eeqnn
Combining this with \eqref{3.4} we know that
 \beqlb\label{3.5}
 \ar\ar
\int_0^t\int_{\mbb{R}}g(s,Y_s(x))\sigma(Y_s(x))^{-\frac12}
1_{\{\sigma(Y_s(x))\neq0\}} M(ds,dx) \cr
 \ar\ar\qquad=
\int_0^t\int_0^\infty g(s,u)1_{\{u\le Y_s(\infty)\}}W(ds,du), \qquad
g\in \mathcal{F}_0.
 \eeqlb

{\bf Step 4.} Observe that for each $x\in\mbb{R}$,
 \beqlb\label{3.6}
 \ar\ar
\int_0^t\int_{\mbb{R}}1_{\{y\le x\}}M(ds,dy) \cr
 \ar=\ar
\int_0^t\int_{\mbb{R}}1_{\{Y_s(y)\le Y_s(x)\}}M(ds,dy)
 -\int_0^t\int_{\mbb{R}}1_{\{Y_s(y)\le Y_s(x),y>x\}}M(ds,dy) \cr
 \ar=\ar
\int_0^t\int_{\mbb{R}}\Big[1_{\{Y_s(y)\le
Y_s(x)\}}\sigma(Y_s(y))^{\frac12}\Big]\sigma(Y_s(y))^{-\frac12}
1_{\{\sigma(Y_s(y))\neq0\}}M(ds,dy)\cr
 \ar\ar
 +
\int_0^t\int_{\mbb{R}}1_{\{Y_s(y)\le
Y_s(x),\sigma(Y_s(y))=0\}}M(ds,dy)
-\int_0^t\int_{\mbb{R}}1_{\{Y_s(y)\le Y_s(x),y>x\}}M(ds,dy) \cr
 \ar=:\ar
\tilde{M}_1(t,x)+\tilde{M}_2(t,x)-\tilde{M}_3(t,x).
 \eeqlb
Letting $g(s,u)=1_{\{u\le Y_s(x)\}}\sigma(u)^\frac12$ in \eqref{3.5}
we obtain
 \beqlb\label{3.7}
\tilde{M}_1(t,x)= \int_0^t\int_0^\infty 1_{\{u\le
Y_s(x)\}}\sigma(u)^\frac12W(ds,du).
 \eeqlb
Observe that
 \beqnn
\mbf{E}\{\tilde{M}_2(t,x)^2\}=
\mbf{E}\Big\{\int_0^tds\int_{\mbb{R}}1_{\{Y_s(y)\le
Y_s(x),\sigma(Y_s(y))=0\}}\sigma(Y_s(y))dy\Big\}=0,
 \eeqnn
which implies that
 \beqlb\label{3.8}
\tilde{M}_2(t,x)=0,\quad\mbf{P}\mbox{-a.s.}
 \eeqlb
Observe that the quadratic variation process $\<\tilde{M}_3(x)\>_t$
of $t\mapsto \tilde{M}_3(t,x)$ satisfies
 \beqnn
\<\tilde{M}_3(x)\>_t
 \ar=\ar
\int_0^tds\int_{\mbb{R}}\sigma(X_s(-\infty,y]))1_{\{Y_s(y)\le
Y_s(x),y>x\}}X_s(dx) \cr
 \ar=\ar
\int_0^tds\int_{\mbb{R}}\sigma(Y_s(y))1_{\{Y_s(y)\le
Y_s(x),y>x\}}dY_s(y)=0,
 \eeqnn
which implies $\tilde{M}_3(t,x)=0$  $\mbf{P}$-a.s. Combining
\eqref{3.6}--\eqref{3.8} one has
 \beqnn
\int_0^t\int_{\mbb{R}}1_{\{y\le x\}}M(ds,dy) =\int_0^t\int_0^\infty
1_{\{u\le Y_s(x)\}}\sigma(u)^\frac12W(ds,du), \quad
\mbf{P}\mbox{-a.s.}
 \eeqnn
Then by stochastic Fubini's theorem, for $f\in C_0^2(\IR)$ with
$\<|f|,1\><\infty$ and $\bar{f}(y):=\int_y^\infty f(x)dx$,
 \beqnn
M_t(\bar{f})
 \ar=\ar
\int_0^t\int_{\mbb{R}}\bar{f}(y)M(ds,dy)
=\int_0^t\int_{\mbb{R}}\Big[\int_{\mbb{R}}f(x)1_{\{y\le
x\}}dx\Big]M(ds,dy) \cr
 \ar=\ar
\int_{\mbb{R}}f(x)dx\int_0^t\int_{\mbb{R}}1_{\{y\le x\}}M(ds,dy)\cr
 \ar=\ar
\int_{\mbb{R}}f(x)dx\int_0^t\int_0^\infty 1_{\{u\le
Y_s(x)\}}\sigma(u)^\frac12W(ds,du), \quad \mbf{P}\mbox{-a.s.}
 \eeqnn
By integration by parts, we have $\mbf{P}$-a.s.
 \beqnn
M_t(\bar{f})=X_t(\bar{f})-X_0(\bar{f})-\frac12\int_0^tX_s((\bar{f})'')ds
=\<Y_t,f\>-\<Y_0,f\>-\frac12\int_0^t\<Y_s,f''\>ds,
 \eeqnn
which yields that $\{Y_t:t\ge0\}$ solves \eqref{3.1}.
 \qed

\subsection{Properties}

In this subsection we discuss the existence of density field and
long time survival-extinction behaviors of the
interacting super-Brownian motion.

\btheorem\label{t3.2}  Suppose that $\{X_t:t\ge0\}$ is the
interacting super-Brownian motion with distribution-function
process $\{Y_t:t\ge0\}$ solving SPDE \eqref{3.1}. Then for any $f\in
C^2(\mbb{R})$, $\mbf{P}$-a.s.
 \beqlb\label{3.12}
X_t(f) = X_0(f)+\frac{1}{2}\int_0^tX_s(f'')ds+\int_0^t\int_0^\infty
f(Y_s^{-1}(u))\sqrt{\sigma(u)}W(ds,du),\quad t\ge0,
 \eeqlb
and for any $f\in B(\mbb{R})$ and $t>0$,
 \beqlb\label{3.13}
X_t(f) = X_0(P_tf) +\int_0^t\int_0^\infty
P_{t-s}f(Y_s^{-1}(u))\sqrt{\sigma(u)}W(ds,du)\quad
\mbf{P}\mbox{-a.s.},
 \eeqlb
where $P_t$ denotes the transition semigroup for Brownian motion and
recall that $Y_t^{-1}(u)$ is defined by \eqref{3.11}. Moreover, for
each $t>0$ the random measure $X_t(dx)$ is absolutely continuous
with respect to the Lebesgue measure and the density $X_t(x)$
satisfies
 \beqlb\label{3.14}
X_t(x) = \int_{\mbb{R}}p_t(x-y)X_0(dy)+\int_0^t\int_0^\infty
p_{t-s}(x-Y_s^{-1}(u))\sqrt{\sigma(u)}W(ds,du)~\mbf{P}\mbox{-a.s.},
 \eeqlb
where $p_t$ is the transition density function for Brownian motion.
 \etheorem
 \proof
By the proof of Lemma \ref{l3.1}, we easily get \eqref{3.12}. We then
finish the proof in the following two steps.

{\bf Step 1.} Let $t>0$ be fixed. For $n\ge1$ define stopping time
$\tau_n$ by
 \beqnn
\tau_n=\inf\{s\ge0:X_s(1)\ge n\}.
 \eeqnn
Then $\tau_n\to\infty$ as $n\to\infty$. In this step we want to show that
for each $f\in C^2(\IR)$ and $n\ge1$,
 \beqlb\label{3.17}
X_{t\wedge \tau_n}(P_{t-(t\wedge\tau_n)}f) =X_0(P_tf)+
\int_0^t\int_0^\infty P_{t-s}f(Y_s^{-1}(u))\sqrt{\sigma(u)}1_{\{s\le
\tau_n\}}W(ds,du) ~\mbf{P}\mbox{-a.s.}
 \eeqlb
Consider a partition $\Delta_m:=\{0=t_0<t_1<\cdots<t_m=t\}$ of
$[0,t]$. Let $|\Delta_m|:=\max_{1\le i\le m}|t_i-t_{i-1}|$. It is
obvious that $\frac{dP_sf}{ds}=\frac12 \Delta P_sf$ for $s\ge0$. It then
follows from \eqref{3.12} that
 \beqlb\label{3.16}
 \ar\ar
X_{t\wedge \tau_n}(P_{t-(t\wedge\tau_n)}f) \cr
 \ar=\ar
X_0(P_tf) + \sum\limits_{i=1}^m X_{t_i\wedge
\tau_n}\big(P_{t-(t_i\wedge\tau_n)}f
-P_{t-(t_{i-1}\wedge\tau_n)}f\big) \cr
 \ar\ar
+ \sum\limits_{i=1}^m \big[X_{t_i\wedge
\tau_n}(P_{t-(t_{i-1}\wedge\tau_n)}f) - X_{t_{i-1}\wedge
\tau_n}(P_{t-(t_{i-1}\wedge\tau_n)}f)\big] \cr
 \ar=\ar
X_0(P_tf) +
\frac12\sum\limits_{i=1}^m\int_{t-(t_{i-1}\wedge\tau_n)}^{t-(t_i\wedge\tau_n)}
X_{t_i\wedge \tau_n}(\Delta P_sf)ds +
\frac12\sum\limits_{i=1}^m\int_{t_{i-1}}^{t_i}
X_s(\Delta P_{t-(t_{i-1}\wedge\tau_n)}f)1_{\{s\le\tau_n\}}ds \cr
 \ar\ar
+ \sum\limits_{i=1}^m\int_{t_{i-1}}^{t_i}\int_0^\infty
P_{t-(t_{i-1}\wedge\tau_n)}f(Y_s^{-1}(u))\sqrt{\sigma(u)}1_{\{s\le\tau_n\}}W(ds,du)
\cr
 \ar=\ar
X_0(P_tf) + \frac12\int_0^t\sum\limits_{i=1}^mI_i(s)
\big[X_s(\Delta P_{t-(t_{i-1}\wedge\tau_n)}f) -X_{t_i\wedge
\tau_n}(\Delta P_{t-s}f)\big]1_{\{s\le\tau_n\}}ds \cr
 \ar\ar
+\int_0^t\int_0^\infty\sum\limits_{i=1}^mI_i(s)
P_{t-(t_{i-1}\wedge\tau_n)}f(Y_s^{-1}(u))\sqrt{\sigma(u)}1_{\{s\le\tau_n\}}W(ds,du)\cr
 \ar=:\ar
X_0(P_tf)+ \tilde{I}_1(t,m,n)/2+\tilde{I}_2(t,m,n),
 \eeqlb
where $I_i(s):=1_{(t_{i-1},t_i]}(s)$. Observe that by the dominated
convergence and the continuities of $s\mapsto P_sf$ and $s\mapsto
X_s(P_{t'}f)$ for $t'>0$, we have
 \beqnn
\lim_{|\Delta_m|\to0}|\tilde{I}_1(t,m,n)|
 \ar\le\ar
\int_0^t\lim_{|\Delta_m|\to0}\sum\limits_{i=1}^mI_i(s)
\Big[\big|X_s(\Delta P_{t-(t_{i-1}\wedge\tau_n)}f) -X_s(\Delta P_{t-s}f)\big|
\cr
 \ar\ar\qquad\qquad
+\big|X_s(\Delta P_{t-s}f) -X_{t_i\wedge
\tau_n}(\Delta P_{t-s}f)\big|\Big]1_{\{s\le\tau_n\}}ds
 \eeqnn
and
 \beqnn
 \ar\ar
\lim_{|\Delta_m|\to0}\mbf{E}\Big\{\Big[\tilde{I}_2(t,m,n)-\int_0^t\int_0^\infty
P_{t-s}f(Y_s^{-1}(u))\sqrt{\sigma(u)}1_{\{s\le\tau_n\}}W(ds,du)\Big]^2\Big\}
\cr
 \ar\ar\qquad=
\lim_{|\Delta_m|\to0}\mbf{E}\Big\{\Big[\int_0^t\int_0^\infty\sum\limits_{i=1}^mI_i(s)
\big[P_{t-(t_{i-1}\wedge\tau_n)}f(Y_s^{-1}(u)) \cr
 \ar\ar\qqquad\qqquad\qquad
-P_{t-s}f(Y_s^{-1}(u))\big]
\sqrt{\sigma(u)}1_{\{s\le\tau_n\}}W(ds,du)\Big]^2\Big\} \cr
 \ar\ar\qquad=
\mbf{E}\Big\{\int_0^tds\int_0^\infty\lim_{|\Delta_m|\to0}\sum\limits_{i=1}^mI_i(s)\big[
P_{t-(t_{i-1}\wedge\tau_n)}f(Y_s^{-1}(u)) \cr
 \ar\ar\qqquad\qqquad\qqquad\qquad
-P_{t-s}f(Y_s^{-1}(u))\big]^2 \sigma(u)1_{\{s\le\tau_n\}}du\Big\}=0.
 \eeqnn
Thus \eqref{3.17} follows by letting $|\Delta_m|\to0$ in
\eqref{3.16}.

{\bf Step 2}. Equation \eqref{3.17} holds for all $f\in B(\IR)$ by
an approximation argument and then \eqref{3.13} follows by letting
$n\to\infty$. We now prove the last assertion. By \eqref{3.17} we
have that for any $f\in B(\IR)$,
 \beqnn
\mbf{E}\Big\{\<X_{t\wedge \tau_n},P_{t-(t\wedge\tau_n)}f\>\Big\}
=X_0(P_tf).
 \eeqnn
It follows from Fatou's lemma that for $f\in B(\IR)^+$,
 \beqnn
\mbf{E}\{X_t(f)\}=\mbf{E}\Big\{\liminf_{n\to\infty}\<X_{t\wedge
\tau_n},P_{t-(t\wedge\tau_n)}f\>\Big\}
 \le
\liminf_{n\to\infty}\mbf{E}\Big\{\<X_{t\wedge
\tau_n},P_{t-(t\wedge\tau_n)}f\>\Big\}=X_0(P_tf).
 \eeqnn
Thus,
 \beqnn
 \ar\ar
\mbf{E}\Big\{\int_{\mbb{R}}f(x)dx\int_0^t1_{\{s\le
\tau_n\}}ds\int_0^\infty p_{t-s}(x-Y_s^{-1}(u))^2\sigma(u)du\Big\}
\cr
 \ar\ar\quad\le
\mbf{E}\Big\{\int_{\mbb{R}}f(x)dx\int_0^t(t-s)^{-\frac12}ds\int_0^\infty
p_{t-s}(x-Y_s^{-1}(u))\sigma(u)1_{\{s\le \tau_n\}}du\Big\} \cr
 \ar\ar\quad\le
\mbf{E}\Big\{\int_{\mbb{R}}f(x)dx\int_0^t(t-s)^{-\frac12}ds\int_0^\infty
p_{t-s}(x-y)\sigma(Y_s(y))1_{\{s\le \tau_n\}}X_s(dy)\Big\} \cr
 \ar\ar\quad\le
c_n\mbf{E}\Big\{\int_{\mbb{R}}f(x)dx\int_0^t(t-s)^{-\frac12}ds\int_0^\infty
p_{t-s}(x-y)X_s(dy)\Big\} \cr
 \ar\ar\quad\le
c_n\mbf{E}\Big\{\int_0^t(t-s)^{-\frac12}X_s(P_{t-s}f)ds\Big\} \le
c_n\mu(P_tf)\int_0^t(t-s)^{-\frac12}ds<\infty,
 \eeqnn
where $c_n:=\sup_{x\in[0,n]}\sigma(x)$. Then by \cite[Theorem
7.24]{Li11}, for each $t_0>0$ we have
 \beqlb\label{3.15}
\<X_{t_0\wedge
\tau_n},P_{t_0-(t_0\wedge\tau_n)}f\>=\int_{\mbb{R}}f(x)K_{t_0}(t_0\wedge\tau_n,x)dx
\quad \mbf{P}\mbox{-a.s.}
 \eeqlb
where
 \beqnn
K_{t_0}(t,x):=\int_{\mbb{R}}p_{t_0}(x-y)X_0(dy)+
\int_0^t\int_0^\infty
p_{t_0-s}(x-Y_s^{-1}(u))\sqrt{\sigma(u)}W(ds,du).
 \eeqnn
Letting $n\to\infty$ in \eqref{3.15} we get
$X_t(f)=\int_{\mbb{R}}f(x)K_t(t,x)dx$ $\mbf{P}$-a.s., which finishes
the proof.
 \qed

Write $\sigma_0(x):=\int_0^x \sigma(y) dy$ for $x\geq 0$. Notice
that $\sigma_0(x)$ is an increasing continuous function and
$\sigma_0(x)=0$ if and only if $x=0$. In the following of this
subsection we always assume that $\mu\in M(\mbb{R})$ and
$\{X_t:t\ge0\}$ is an interacting super-Brownian motion with
$X_0=\mu$. Taking $f=1$ in \eqref{3.12} we get
 \beqnn
X_t(1)=X_0(1)+\int_0^t\int_0^{X_s(1)} \sqrt{\sigma(u)}W(ds,du).
 \eeqnn
Thus the total mass process $\{X_t(1):t\ge0\}$ is a nonnegative
continuous local martingale with quadratic variation process
 \beqlb\label{3.9}
\< X(1)\>_t = \int_0^tds\int_0^{X_s(1)}\sigma(y)dy=\int_0^t
\sigma_0(X_s(1)) d s.
 \eeqlb
By the uniqueness of solution to the local martingale problem, if
$\mu$ is a zero measure, then $\mbf{P}\{X_t(1)=0 \,\,\text{for
all}\,\, t>0\}=1$. So we assume that $\mu(1)>0$ in the
rest of the subsection.

For $a\geq 0$ let
 \beqnn
\hat{\tau}_a:=\inf\{t\geq 0: X_t(1)=a\}
 \eeqnn
with the convention $\inf\emptyset=\infty$.

\btheorem\label{t3.4} For any finite measure $\mu$  and constant
$a>0$ satisfying $\mu(1)>a>0$, we have $\mbf{P}\{\hat{\tau}_a<
\infty\}=1$. Further, $\mbf{P}$-a.s. $X_t(1)\rightarrow 0$ as
$t\rightarrow\infty$. \etheorem

\begin{proof}
Notice that for any $\tilde{\la}>0$,  process
 \beqnn \exp\Big\{-\tilde{\la}
X_t(1)-\frac{\tilde{\la}^2}{2}\int_0^t \sigma_0(X_s(1)) d s\Big\}
 \eeqnn
 is a bounded continuous martingale. Then by optional sampling,
\beqnn e^{-\tilde{\la} \mu(1)}=\mbf{E}\Big\{ \exp\Big[-\tilde{\la}
X_{t\wedge\hat{\tau}_a}(1)
-\frac{\tilde{\la}^2}{2}\int_0^{t\wedge\hat{\tau}_a}\sigma_0(X_s(1))
d s\Big]\Big\}. \eeqnn
 Letting $t\rightarrow\infty$, by dominated convergence we have
 \beqnn
 \ar\ar
\mbf{E}\Big\{ \exp\Big[-\tilde{\la}
a-\frac{\tilde{\la}^2}{2}\int_0^{\hat{\tau}_a} \sigma_0(X_s(1)) d
s\Big]\Big\} \cr
 \ar\ar\quad=
\lim_{t\goto\infty} \mbf{E}\Big\{ \exp\Big[-\tilde{\la}
X_{t\wedge\hat{\tau}_a}(1)-\frac{\tilde{\la}^2}{2}\int_0^{t\wedge\hat{\tau}_a}\sigma_0(
X_s(1)) d s\Big]\Big\} = e^{-\tilde{\la} \mu(1)}.
 \eeqnn
Now letting $\tilde{\la}\rightarrow 0$ we have
 \beqnn
\mbf{P}\Big\{\int_0^{\hat{\tau}_a} \sigma_0( X_s(1)) d
s<\infty\Big\}= 1.
 \eeqnn
Observe that by the assumption on function $\sigma$ we have $\inf_{y\geq x}\sigma_0(y)>0$ for any $x>0$,  which
implies that
 \beqnn
\{\hat{\tau}_a=\infty\}\subset \Big\{\int_0^{\hat{\tau}_a}
\sigma_0(X_s(1)) d s=\infty\Big\}.
 \eeqnn
Then we have $\mbf{P}\{\hat{\tau}_a<\infty\}=1$, which proves the
first assertion.

For any $n\ge1$ and small $\ep>0$ let
 \beqnn
A_n:=\{\hat{\tau}_{\ep^{n+1}}<\infty, \,\, \hat{\tau}_\ep\circ
\theta_{\hat{\tau}_{\ep^{n+1}}} <\hat{\tau}_{\ep^{n+2}}\circ
\theta_{\hat{\tau}_{\ep^{n+1}}} \},
 \eeqnn
where $\theta_t, \, t\ge 0,$ denotes the usual shift operator,
that is, $Y\equiv\theta_t(X)$ is the process  such that $Y_s=X_{t+s}$ for all $s\ge0$.

Since $\{X_t(1):t\ge0\}$ is a continuous local martingale, then by
optional stopping we have
 \beqnn
\ep^{n+1}
 \ar=\ar
X_{\hat{\tau}_{\ep^{n+1}}}(1)=\mbf{E}\big[X_{\hat{\tau}_{1,\ep,n}
\wedge\hat{\tau}_{2,\ep,n}}(1)
\big|\mathscr{F}_{\hat{\tau}_{\ep^{n+1}}}\big] \cr
 \ar\ge\ar
\mbf{E}\big[X_{\hat{\tau}_{1,\ep,n}}(1)1_{\{\hat{\tau}_{1,\ep,n}<\hat{\tau}_{2,\ep,n}\}}
\big|\mathscr{F}_{\hat{\tau}_{\ep^{n+1}}}\big] =
\ep\mbf{P}\big[\hat{\tau}_{1,\ep,n}<\hat{\tau}_{2,\ep,n}
\big|\mathscr{F}_{\hat{\tau}_{\ep^{n+1}}}\big],
 \eeqnn
where $\hat{\tau}_{1,\ep,n}:=\hat{\tau}_\ep\circ
\theta_{\hat{\tau}_{\ep^{n+1}}}$ and
$\hat{\tau}_{2,\ep,n}:=\hat{\tau}_{\ep^{n+2}}\circ
\theta_{\hat{\tau}_{\ep^{n+1}}}$. Then we know that
 \beqnn
\mbf{P}\{A_n\}\le\mbf{E}\Big\{\mbf{P}\big[\hat{\tau}_{1,\ep,n}<\hat{\tau}_{2,\ep,n}
\big|\mathscr{F}_{\hat{\tau}_{\ep^{n+1}}}\big]\Big\} \leq \ep^n.
 \eeqnn
It follows from the Borel-Cantelli lemma that $\mbf{P}\{A_n \,\,
\text{infinitely often}\}=0$. Then by the first assertion, for all
$n$ large enough,
 \beqnn
\hat{\tau}_{\ep^{n}}<\infty,\quad\hat{\tau}_{\ep^{n+1}}\circ
\theta_{\hat{\tau}_{\ep^{n}}}<\hat{\tau}_\ep\circ
\theta_{\hat{\tau}_{\ep^{n}}},\qquad \mbf{P}\mbox{-a.s.}
 \eeqnn
Therefore, $\mbf{P}$-a.s. for $n$ large enough $X_t(1)<\ep$  for all
$t\in [\hat{\tau}_{\ep^n}, \hat{\tau}_0)$, where
$\hat{\tau}_0=\lim_{m\to\infty}\hat{\tau}_{\ep^m}$. If
$\hat{\tau}_0<\infty$, then $X_t(1)=0$ for all $t\geq\hat{\tau}_0$;
otherwise,
 $X_t(1)\goto 0$ as $t\goto \infty$ but $X_t(1)>0$ for all $t$.
 Putting them together, we have $X_t(1)\goto 0$.
 \qed
\end{proof}

The extinction behavior of $\{X_t:t\ge0\}$ depends on the branching
rate when the total mass is close to $0$. So, it depends on how fast
$\sigma_0(x)$ converges to $0$ when $x\goto 0+$.

\btheorem\label{t3.6} (i) If there exists a constant
$\gamma_1\in[0,2)$ so that $\liminf_{x\goto 0+}
\sigma_0(x)/x^{\gamma_1}>0$, we have
$\mbf{P}\{\hat{\tau}_0<\infty\}=1$ and $X_t(1)=0$ for all large
enough $t$. (ii) If there exists a constant $\gamma_2\in[2,\infty)$
so that $\limsup_{x\goto 0+} \sigma_0(x)/x^{\gamma_2}<\infty$, we
have $\mbf{P}\{\hat{\tau}_0=\infty\}=1$ and $X_t(1)>0$ for all
$t>0$.
 \etheorem

\bcorollary If there is a constant $\gamma'\ge0$ so that
$\sigma(x)=x^{\gamma'}$ for all $x\geq 0$, then
$\mbf{P}\{\hat{\tau}_0<\infty\}=1$ for $\gamma'<1$ and
$\mbf{P}\{\hat{\tau}_0=\infty\}=1$ for $\gamma'\geq 1$. \ecorollary
 \proof
Since $\sigma_0(x)=\int_0^x \sigma(y)dy=x^{\gamma'+1}/(\gamma'+1)$,
the assertion follows from Theorem \ref{t3.6} immediately.
 \qed

We end this section with the proof of Theorem \ref{t3.6}.

\noindent{\it Proof of Theorem \ref{t3.6}.}
(i) Fix an $\ep>0$ small enough so that $\sigma_0(x)\geq b_1
x^{\gamma_1}$ for all $ 0\leq x\leq \sqrt{\ep}$ and some $b_1>0$,
let $T':=\hat{\tau}_{\ep^{1+\de}}\wedge \hat{\tau}_{\ep^{1/2}}$ for
$\delta>0$. By \eqref{3.9} one can see that
 \beqnn
\exp\Big\{-\tilde{\la}
X_t(1)-\frac{\tilde{\la}^2}{2}\<X(1)\>_t\Big\} =
\exp\Big\{-\tilde{\la} X_t(1)-\frac{\tilde{\la}^2}{2}\int_0^t
\sigma_0(X_s(1))  d s\Big\}
 \eeqnn
is a continuous local martingale. Then by optional sampling, for
$\mu(1)=\ep$ we have
 \beqlb\label{3.10}
e^{-\tilde{\la}\ep}
 \ar=\ar
\mbf{E}\Big\{ \exp\Big[-\tilde{\la}
X_{T'}(1)-\frac{\tilde{\la}^2}{2}\int_0^{T'} \sigma_0(X_s(1))  d
s\Big]\Big\} \cr
 \ar\le\ar
\mbf{E}\Big\{ \exp\Big[-\tilde{\la}
\ep^{1+\de}-\frac{\tilde{\la}^2}{2}b_1\ep^{\gamma_1(1+\de)}\hat{\tau}_{\ep^{1+\de}}\Big]
1_{\{\hat{\tau}_{\ep^{1+\de}}<\hat{\tau}_{\ep^{1/2}}\}}\Big\}
+\mbf{P}\{\hat{\tau}_{\ep^{1+\de}}>\hat{\tau}_{\ep^{1/2}}\}.
 \eeqlb
Let
$\hat{\tau}_{3,\ep}:=\hat{\tau}_{\ep^{1+\de}}\wedge\hat{\tau}_{\ep^{1/2}}$.
Then by optional stopping again,
 \beqnn
\ep=\mbf{E}\{X_{\hat{\tau}_{3,\ep}}(1)\}\ge\mbf{E}\Big\{X_{\hat{\tau}_{\ep^{1/2}}}(1)
1_{\{\hat{\tau}_{\ep^{1+\de}}>\hat{\tau}_{\ep^{1/2}}\}}\Big\}
=\sqrt{\ep}\mbf{P}\{\hat{\tau}_{\ep^{1+\de}}>\hat{\tau}_{\ep^{1/2}}\},
 \eeqnn
which implies
 \beqlb\label{extinctA}
\mbf{P}\{\hat{\tau}_{\ep^{1+\de}}>\hat{\tau}_{\ep^{1/2}}\} \le
\sqrt{\ep}.
 \eeqlb

Since $0\leq \gamma_1 <2 $, we can choose $0<\de<1/2$ so that
\begin{equation}\label{extinctB}
\gamma_1(1+\de)<2-4\de.
\end{equation}
Then by \eqref{3.10}--\eqref{extinctB}, for
$\tilde{\la}=\ep^{-1+\de}$ we have
 \beqlb\label{extinctC}
1-2\ep^\de
 \ar\leq\ar  e^{-\tilde{\la}(\ep-\ep^{1+\de})}\cr
 \ar\leq\ar
\mbf{E}\Big\{
\exp\big[-2^{-1}\tilde{\la}^2b_1\ep^{\gamma_1(1+\de)}\hat{\tau}_{\ep^{1+\de}}\big]
1_{\{\hat{\tau}_{\ep^{1+\de}}<\hat{\tau}_{\ep^{1/2}}\}}\Big\}
+\mbf{P}\{\hat{\tau}_{\ep^{1+\de}}>\hat{\tau}_{\ep^{1/2}}\}e^{\tilde{\la}
\ep^{1+\de}} \cr
 \ar\le\ar
\mbf{E}\Big\{ \exp\big[-
2^{-1}b_1\ep^{-2\de}\hat{\tau}_{\ep^{1+\de}}\big]\Big\}+\sqrt{\ep}
e^{\tilde{\la} \ep^{1+\de}} \cr
 \ar\leq\ar
\mbf{P}\{\hat{\tau}_{\ep^{1+\de}} \leq\ep^\de
\}+\exp\big[-2^{-1}b_1\ep^{-2\de}\ep^\de\big]
\mbf{P}\{\hat{\tau}_{\ep^{1+\de}}> \ep^\de
\}+\sqrt{\ep}e^{\tilde{\la} \ep^{1+\de}} \cr
 \ar=\ar
1-\mbf{P}\{\hat{\tau}_{\ep^{1+\de}}> \ep^\de
\}+\exp\big[-2^{-1}b_1\ep^{-\de}\big]\mbf{P}\{\hat{\tau}_{\ep^{1+\de}}>
\ep^\de \}+\sqrt{\ep}e^{\tilde{\la} \ep^{1+\de}}.
 \eeqlb
Solving inequality (\ref{extinctC}) for
$\mbf{P}\{\hat{\tau}_{\ep^{1+\de}}> \ep^\de \}$, then for small
$\ep$,
 \beqnn
\mbf{P}\{\hat{\tau}_{\ep^{1+\de}}> \ep^\de
\}\leq\frac{2\ep^\de+\sqrt{\ep}e^{\tilde{\la}
\ep^{1+\de}}}{1-\exp\big[-2^{-1}b_1\ep^{-\de}\big]}\leq 3\ep^\de.
 \eeqnn

Put $x_n:=\ep^{(1+\de)^n}$ for $n\ge1$. Then by the Markov property
of $\{X_t(1):t\ge0\}$,
 \beqnn
\mbf{P}\{\hat{\tau}_{x_{n+1}}-\hat{\tau}_{x_n}>
x_n^\de|\mathscr{F}_{\hat{\tau}_{x_n}} \}
 =
\mbf{P}\{\hat{\tau}_{x_{n+1}}\circ \theta_{\hat{\tau}_{x_n}}>
x_n^\de|\mathscr{F}_{\hat{\tau}_{x_n}} \}\leq 3x_n^\de,
 \eeqnn
which implies that
 \beqnn
\sum_{n=0}^\infty\mbf{P}\{\hat{\tau}_{x_{n+1}}-\hat{\tau}_{x_n}>
x_n^\de\} =\sum_{n=0}^\infty
\mbf{E}\Big\{\mbf{P}\big[\hat{\tau}_{x_{n+1}}-\hat{\tau}_{x_n}>
x_n^\de|\mathscr{F}_{\hat{\tau}_{x_n}}\big]\Big\} \leq
\sum_{n=0}^\infty 3x_n^\de<\infty.
 \eeqnn
Then by the Borel-Cantelli lemma we have
 \beqnn
\mbf{P}\Big\{ \{\hat{\tau}_{x_{n}}<\infty\}\cap\{
\hat{\tau}_{x_{n+1}}-\hat{\tau}_{x_{n}}
>x_n^\de\} \,\,\, \text{infinitely often}\Big\}=0.
 \eeqnn
It follows from Theorem \ref{t3.4} that there are at most finitely
many $n$ so that $\hat{\tau}_{x_n}<\infty$ and
$\hat{\tau}_{x_{n+1}}-\hat{\tau}_{x_n}>x_n^\de$, $\mbf{P}$-a.s.
Since $\mbf{P}\{ \hat{\tau}_{x_{n}}<\infty$ for all $n\}=1$ and
$\sum_{n=0}^\infty x_n^\de<\infty$, then
 \beqnn
\hat{\tau}_0=\lim_{n\to\infty}\hat{\tau}_{x_n}
=\hat{\tau}_{x_0}+\sum_{n=0}^\infty[\hat{\tau}_{x_{n+1}}-\hat{\tau}_{x_{n}}]
\le \hat{\tau}_{x_0}+\sum_{n=0}^\infty x_n^\delta <\infty, \quad
\mbf{P}\mbox{-a.s.}
 \eeqnn
and we have the desired result for $\mu $ with $\mu(1)=\ep$. Then
the result for general $\mu$ follows.

(ii) Given small $\ep, \delta>0$ so that $\sigma_0(x)\leq b_2
x^{\gamma_2}$ for all $0<x\leq \ep^{1-\de} $ and some $b_2>0$,
 for any finite measure $\mu$ on $\bR$
with $\mu(1)=\ep$ and $T'':=\hat{\tau}_{\ep^{1+\de}} \wedge
\hat{\tau}_{\ep^{1-\delta}}$. By Ito's formula we have
 \beqnn
X_{t\wedge T''}(1)^{-1} =X_0(1)^{-1}-\int_0^{t\wedge T''}
X_s(1)^{-2} d X_s(1) +\int_0^{t\wedge T''}
\sigma_0(X_s(1))X_s(1)^{-3} d s.
 \eeqnn
Then by integration by parts,
 \beqnn
 \ar\ar
X_{t\wedge T''}(1)^{-1}\exp\Big\{-\int_0^{t\wedge
T''}\sigma_0(X_s(1)) X_s(1)^{-2} d s\Big\}\cr
 \ar\ar\quad
=X_0(1)^{-1}+\int_0^t\exp\Big\{-\int_0^{s\wedge T''}\sigma_0(X_r(1))
X_r(1)^{-2}d r\Big\} d (X_{s\wedge T''}(1)^{-1}) \cr
 \ar\ar\quad\quad
+\int_0^t  X_{s\wedge T''}(1)^{-1} d \Big(\exp\Big\{-\int_0^{t\wedge
T''}\sigma_0(X_s(1))X_s(1)^{-2} d s\Big\} \Big)\cr
 \ar\ar\quad
=X_0(1)^{-1}-\int_0^{t\wedge T''}X_s(1)^{-2}
\exp\Big\{-\int_0^{t\wedge T''}\sigma_0(X_s(1))X_s(1)^{-2} d s\Big\}
d X_s(1)
 \eeqnn
is a martingale. Then by optional sampling
 \beqnn
\ep^{-1}
 \ar=\ar
\mbf{E}\Big\{X_{t\wedge T''}(1)^{-1} \exp\Big[ -\int_0^{t\wedge T''}
\sigma_0(X_s(1))X_s(1)^{-2}d s\Big]\Big\}\cr
 \ar\geq\ar
\mbf{E}\Big\{X_{t\wedge \hat{\tau}_{\ep^{1+\de}}}(1)^{-1}
\exp\Big[-\int_0^{t\wedge \hat{\tau}_{\ep^{1+\de}}}
b_2X_{s}(1)^{\gamma_2-2} d
s\Big]1_{\{\hat{\tau}_{\ep^{1+\de}}<\hat{\tau}_{\ep^{1-\de}}\}}
\Big\}\cr
 \ar\geq\ar
\mbf{E}\Big\{X_{t\wedge \hat{\tau}_{\ep^{1+\de}}}(1)^{-1}\exp\big[-
b_2\ep^{(1-\de)(\gamma_2-2)}\hat{\tau}_{\ep^{1+\de}} \big]
1_{\{\hat{\tau}_{\ep^{1+\de}}<\hat{\tau}_{\ep^{1-\de}}\}}\Big\}.
 \eeqnn
Letting $t\to\infty$ in above inequality we get
 \beqnn
\ep^{-1}
 \ar\geq\ar
\mbf{E}\Big\{X_{\hat{\tau}_{\ep^{1+\de}}}(1)^{-1}\exp\big[-
b_2\ep^{(1-\de)(\gamma_2-2)}\hat{\tau}_{\ep^{1+\de}} \big]
1_{\{\hat{\tau}_{\ep^{1+\de}}<\hat{\tau}_{\ep^{1-\de}}\}}\Big\} \cr
  \ar=\ar
\mbf{E}\Big\{\ep^{-1-\de}\exp\big[-
b_2\ep^{(1-\de)(\gamma_2-2)}\hat{\tau}_{\ep^{1+\de}} \big]
1_{\{\hat{\tau}_{\ep^{1+\de}}<\hat{\tau}_{\ep^{1-\de}}\}}\Big\},
 \eeqnn
which implies
\begin{equation}\label{est_a}
\mbf{E}\Big\{\exp\big[-
b_2\ep^{(1-\de)(\gamma_2-2)}\hat{\tau}_{\ep^{1+\de}} \big]
1_{\{\hat{\tau}_{\ep^{1+\de}}<\hat{\tau}_{\ep^{1-\de}}\}}\Big\}
 \leq
\ep^\de.
\end{equation}

Applying optional sampling we have
 \beqnn
\ep=\mbf{E}\{X_{T''}(1)\}\ge
\mbf{E}\Big\{X_{\hat{\tau}_{\ep^{1-\de}}}(1)
1_{\{\hat{\tau}_{\ep^{1-\de}}<\hat{\tau}_{\ep^{1+\de}}\}}\Big\}=
\ep^{1-\de}\mbf{P}\{\hat{\tau}_{\ep^{1-\de}}<\hat{\tau}_{\ep^{1+\de}}\},
 \eeqnn
which implies
$\mbf{P}\{\hat{\tau}_{\ep^{1-\de}}<\hat{\tau}_{\ep^{1+\de}}\}\leq
\ep^\de$. It then follows from the Markov inequality that
 \beqnn
 \ar\ar
\mbf{P}\{T''<1\} \cr
 \ar\leq\ar
\mbf{P}\{\hat{\tau}_{\ep^{1-\de}}<\hat{\tau}_{\ep^{1+\de}}\}
+\mbf{P}\{
\hat{\tau}_{\ep^{1+\de}}<1\wedge\hat{\tau}_{\ep^{1-\de}}\} \cr
 \ar\leq\ar
\mbf{P}\{\hat{\tau}_{\ep^{1-\de}}<\hat{\tau}_{\ep^{1+\de}}\}
+\mbf{P}\Big\{\exp\big[-b_2\ep^{(1-\de)(\gamma_2-2)}
{\hat{\tau}_{\ep^{1+\de}}} \big] \cr
 \ar\ar\qquad\qqquad\qqquad
>\exp\big[-b_2\ep^{(1-\de)(\gamma_2-2)} \big],
\hat{\tau}_{\ep^{1+\de}}<1\wedge\hat{\tau}_{\ep^{1-\de}} \Big\} \cr
 \ar\leq\ar
\ep^\de+\exp\big[b_2\ep^{(1-\de)(\gamma_2-2)}\big]
\mbf{E}\Big\{\exp\big[-b_2\ep^{(1-\de)(\gamma_2-2)}
{\hat{\tau}_{\ep^{1+\de}}}\big]
1_{\{\hat{\tau}_{\ep^{1+\de}}<1\wedge\hat{\tau}_{\ep^{1-\de}}\}}\Big\}
\cr
 \ar\leq\ar
\ep^\de+\exp\big[b_2\ep^{(1-\de)(\gamma_2-2)}\big]\ep^\de
 \le
(1+e^{b_2})\ep^\de,
 \eeqnn
where we have used (\ref{est_a}) for the fourth inequity.

Taking $x_n=\ep^{(1+\de)^n}$ and repeating the above argument with
$T''$ replaced by $T_n:=\hat{\tau}_{x_n^{1+\de}}\wedge
\hat{\tau}_{x_n^{1-\de}}$ for $n\ge1$, we have
 \beqnn
 \ar\ar
\mbf{P}\big\{\hat{\tau}_{x_{n+1}}-\hat{\tau}_{x_n}<1|\mathscr{F}_{\hat{\tau}_{x_n}}\big\}
= \mbf{P}\big\{\hat{\tau}_{x_{n+1}}\circ
\theta_{\hat{\tau}_{x_n}}<1|\mathscr{F}_{\hat{\tau}_{x_n}}\big\} \cr
 \ar\ar\qquad\le
\mbf{P}\big\{T_n\circ
\theta_{\hat{\tau}_{x_n}}<1|\mathscr{F}_{\hat{\tau}_{x_n}}\big\}
\leq(1+e^{b_2}) x_n^\de =(1+e^{b_2})\ep^{\de(1+\de)^n}.
 \eeqnn
It then follows that
 \beqnn
 \ar\ar
\mbf{P}\big\{\hat{\tau}_{x_{n}}<\infty, \hat{\tau}_{x_{n+1}}-
\hat{\tau}_{x_{n}}<1\big\}
 \le
\mbf{P}\big\{ \hat{\tau}_{{x_{n+1}}}- \hat{\tau}_{{x_{n}}}<1\big\}
\cr
 \ar\ar\qquad\leq
\mbf{E}\Big\{ \mbf{P}\big[\hat{\tau}_{{x_{n+1}}}-
\hat{\tau}_{{x_{n}}}<1|\mathscr{F}_{\hat{\tau}_{x_n}}\big]\Big\}
 \le
(1+e^{b_2})\ep^{\de(1+\de)^n}.
 \eeqnn
By the Borel-Cantelli lemma,
 \beqnn
\mbf{P}\Big\{\{\hat{\tau}_{x_{n}}<\infty\}\cap\{\hat{\tau}_{x_{n+1}}-
\hat{\tau}_{x_{n}}<1\} \quad \text{infinitely often}\Big\}=0.
 \eeqnn
It thus follows from Theorem \ref{t3.4} that for $n$ large enough,
$\hat{\tau}_{x_{n+1}}- \hat{\tau}_{x_{n}}\geq 1 $. We can thus
conclude that $\mbf{P}\{\hat{\tau}_0=\infty\}=1$. Then the desired
result follows for any positive measure $\mu$. The assertion on
behavior of the total mass process follows immediately from the end
of the proof of Theorem \ref{t3.6}.
 \qed

\section{Interacting Fleming-Viot processes}

\setcounter{equation}{0}

Let $D_1(\mbb{R})$ be the subset of $D(\mbb{R})$ consisting of
continuous functions $f$ with $f(\infty)=1$ and $M_1(\IR)$ be the
space of probability measures on $\mbb{R}$ equipped with the
topology of weak convergence. Then there is an obvious one-to-one
correspondence between $D_1(\mbb{R})$ and $M_1(\mbb{R})$ by
associating a probability measure to its distribution function. We
endow $D_1(\mbb{R})$ with the topology induced by this
correspondence from the weak convergence topology on $M_1(\mbb{R})$.

In this section we study the continuous $M_1(\IR)$-valued solution
to the following martingale problem for \textit{interacting
Fleming-Viot process} $\{X_t:t\ge0\}$: for any $f\in
C^2(\mathbb{R})$,
 \beqlb\label{MP2a}
N_t(f):=X_t(f)-X_0(f)-\frac{1}{2}\int_0^t X_s(f'')ds,\quad t>0
 \eeqlb
is a continuous martingale with quadratic variation process
 \beqlb\label{MP2b}
\<N(f)\>_t=\int_0^tds\int_\mathbb{R}X_s(dx)\int_\mathbb{R}(f(y)-f(x))^2
\gamma_0(X_s,x,y)X_s(dy),
 \eeqlb
where $\gamma_0(\mu,x,y):=\gamma(\mu(-\infty,x],\mu(-\infty,y])$ for
$\gamma\in B([0,1]^2)^+$, $\mu\in M_1(\IR)$ and $x,y\in\IR$. If
$\gamma$ is a positive constant function, then the solution
$\{X_t:t\ge0\}$ is the so-called Fleming-Viot process with constant
resampling rate $\gamma$ and Brownian mutation on type space
$\mathbb{R}$; see \cite{FV79} for details. In the following two
subsections we show that the martingale problem
(\ref{MP2a})--(\ref{MP2b}) is well-posed and investigate some
properties of this process.

\subsection{Well-posedness of the martingale  problem}

The following lemma is on the connection between the martingale problem
(\ref{MP2a})--(\ref{MP2b}) and the distribution-function-valued SPDE.

\blemma\label{t4.2} A continuous $D_1(\mbb{R})$-valued process
$\{Y_t:t\ge0\}$ is the distribution-function process of the
interacting Fleming-Viot process if and only if there is, on an
enlarged probability space, a Gaussian white noise
$\{W(ds,da,db):s\ge0,~a,b\in[0,1]\}$ with intensity $dsdadb$ so that
$\{Y_t:t\ge0\}$ solves equation
 \beqlb\label{4.1}
Y_t(y)= Y_0(y)+\frac{1}{2}\int_0^t\Delta Y_s(y)ds
+\int_0^t\int_0^1\int_0^11_{\{a\leq Y_s(y)\leq
b\}}\sqrt{\gamma(a,b)} W(ds,da,db).
 \eeqlb
 \elemma

\btheorem\label{t3.3}   The martingale problem
(\ref{MP2a})--(\ref{MP2b}) is well posed. \etheorem

\noindent{\it Proof.} Let $Y_0\in D_1(\IR)$. Take $E=[0,1]^2$,
$\pi=$Lebesgue measure on $[0,1]^2$ and
 \beqnn
G(a,b,x)=1_{\{a\leq x\leq b\}}\sqrt{\gamma(a,b)},\qquad
x,a,b\in[0,1].
 \eeqnn
Then the conditions in Theorem \ref{p2.3} hold. In fact, it
is elementary to check that for each $x_1,x_2\in[0,1]$ with $x_1\le
x_2$,
 \beqnn
\int_0^1da\int_0^1G(a,b,x_1)^2db\le\|\gamma\|\int_0^{x_1}da\int_{x_1}^1db
\le\|\gamma\|x_1
 \eeqnn
and
 \beqnn
 \ar\ar
\int_0^1da\int_0^1|G(a,b,x_1)-G(a,b,x_2)|^2db \le
\|\gamma\|\int_0^1da\int_0^1\big|1_{\{a\leq x_1\leq b\}}-1_{\{a\leq
x_2\leq b\}}\big|db \cr
 \ar\ar\qquad\le
\|\gamma\|\int_0^{x_1}da\int_{x_1}^{x_2}db+\|\gamma\|\int_{x_1}^{x_2}da\int^1_{x_2}db
 \le\|\gamma\|[x_2-x_1].
 \eeqnn
This means that conditions \eqref{2.1} and \eqref{2.2} hold. It is
obvious that $G(a,b,0)=0$ for all $a\in(0,1]$ and $b\in[0,1]$, which
implies that condition \eqref{2.17} holds. Therefore, equation
\eqref{4.1} has a pathwise unique continuous $D(\mbb{R})$-valued
solution $\{Y_t:t\ge0\}$ with $Y_0\in D_1(\IR)$.

Observe that the process $\{\bar{Y}_t:t\ge0\}$, defined by
$\bar{Y}_t(x)=1$ for all $x\in\mbb{R}$ and $t\ge0$, is a solution to
(\ref{4.1}) with $\bar{Y}_0=1$. Since $Y_0\le 1$, then by
Proposition \ref{l2.1}(ii), $Y_t\leq \bar{Y}_t=1$ $\mbf{P}$-a.s.,
especially $Y_t(\infty)\leq1$. Moreover, we known from (\ref{2.18})
that
 \beqnn
\mbf{E}\{Y_t(\infty)\}=Y_0(\infty),
 \eeqnn
thus $Y_0(\infty)=1$ implies that $Y_t(\infty)=1$ $\mbf{P}$-a.s. for
each $t>0$, which means that \eqref{4.1} has a unique
strong continuous $D_1(\mbb{R})$-valued solution as $Y_0\in
D_1(\mbb{R})$. Then the assertion follows from Lemma \ref{t4.2}
immediately. \qed

\noindent{\it Proof of Lemma \ref{t4.2}.} The proof is carried out in two steps.

{\bf Step 1.} Suppose that $\{Y_t:t\ge0\}$ is a solution to
\eqref{4.1}, and $\{X_t:t\ge0\}$ is the corresponding measure-valued
process. By integration by parts, for any $f\in C_0^3(\mathbb{R})$,
we have $X_t(f)=-\<Y_t,f'\>$, and so \eqref{4.1} yields
 \beqnn X_t(f)
 \ar=\ar
-\<Y_t,f'\>-\frac{1}{2}\int_0^t\<Y_s,f'''\>ds \cr
 \ar\ar
-\int_0^t\int_0^1\int_0^1\Big[\int_\mathbb{R}1_{\{a\leq Y_s(y)\leq
b\}}f'(y)dy\Big]\sqrt{\gamma(a,b)}W(ds,da,db) \cr
 \ar=\ar
X_0(f)+\frac{1}{2}\int_0^tX_s(f'')ds \cr
 \ar\ar
+\int_0^t\int_0^1\int_0^1\Big[\int_\mathbb{R}1_{\{a\leq Y_s(y)\leq
b\}}f'(y)dy\Big]\sqrt{\gamma(a,b)}W(ds,da,db),
 \eeqnn
which implies that
 \beqnn
\tilde{N}_t(f)
 \ar:=\ar
X_t(f)-X_0(f)-\frac{1}{2}\int_0^tX_s(f'')ds \cr
 \ar=\ar
\int_0^t\int_0^1\int_0^1(f(Y_s^{-1}(b))-f(Y_s^{-1}(a)))\sqrt{\gamma(a,b)}W(ds,da,db)
 \eeqnn
is a square-integrable continuous martingale with
 \beqnn
\<\tilde{N}(f)\>_t
 =
\int_0^tds\int_\mbb{R}X_s(dx)\int_\mbb{R}(f(y)-f(x))^2
\gamma_0(X_s,x,y)X_s(dy),
 \eeqnn
where $Y_t^{-1}(u)$ is defined by \eqref{3.11}. By an approximation
argument, one can see the above relation remains true for any $f\in
C^2(\IR)$. Thus, $\{X_t:t\ge0\}$ is an interacting Fleming-Viot
process.

 {\bf Step 2.} Suppose that
$\{X_t:t\ge0\}$ is an interacting Fleming-Viot process and
$Y_t(x):=X_t(-\infty,x]$ for $x\in\IR$ and $t\ge0$. Let
$\mathcal{S}(\mathbb{R})$ be the Schwartz space and
$\mathcal{S}'(\mathbb{R})$ the space of Schwartz distributions. Let
$f\in \mathcal{S}(\mathbb{R})$ and $\bar{f}(y)=\int_y^\infty
f(x)dx$. Then
 \beqnn
\<Y_t,f\>
 \ar=\ar
X_t(\bar{f})
=X_0(\bar{f})+\frac{1}{2}\int_0^tX_s((\bar{f})'')ds+M_t(\bar{f}) \cr
 \ar=\ar
\<Y_0,f\>+ \frac{1}{2}\int_0^t\<Y_s,f''\>ds+\hat{I}_t(\bar{f}),
 \eeqnn
where $\{\hat{I}_t(\bar{f}):t\ge0\}$ ia a martingale with quadratic
variation process
 \beqnn
 \<\hat{I}(\bar{f})\>_t
=\int_0^tds\int_\mathbb{R}X_s(dx)\int_\mathbb{R}
(\bar{f}(y)-\bar{f}(x))^2\gamma_0(X_s,x,y)X_s(dy),
\qquad t\geq0.
 \eeqnn
Define the $\mathcal{S}'(\mathbb{R})$-valued process $\{N_t:t\ge0\}$
by $N_t(f)=\hat{I}_t(\bar{f})$ for any $f\in
\mathcal{S}(\mathbb{R})$. Then $\{N_t:t\ge0\}$ is a
square-integrable continuous $\mathcal{S}'(\mathbb{R})$-valued
martingale with
 \beqnn
\<N(f)\>_t
 \ar=\ar
\int_0^tds\int_\mathbb{R}\int_\mathbb{R}[\bar{f}(y)-\bar{f}(x)]^2
\gamma_0(X_s,x,y)X_s(dy)X_s(dx)
\cr
 \ar=\ar
\int_0^tds\int_0^1da\int_0^1[\bar{f}(Y_s^{-1}(b))-\bar{f}(Y_s^{-1}(a))]^2
\gamma(a,b)db
\cr
 \ar=\ar
\int_0^tds\int_0^1\int_0^1da\Big[\int_\mathbb{R}1_{\{a\leq
Y_s(y)\leq b\}}f(y)dy\Big]^2\gamma(a,b)db.
 \eeqnn
Let $f\in L^2(\IR)$ be the subset of $\mathscr{B}(\IR)$ consisting
of functions $f$ with $\<f^2,1\><\infty$ and $L^2([0,1]^2)$ be
defined similarly. Similar to the martingale representation theorem
(see \cite[Theorem 3.3.6]{KX95} or \cite[Theorem III-6]{EM90}),
there exist a $L^2([0,1]^2)$-cylindrical Brownian motion
$\{B_t:t\ge0\}$ (possibly on an extended probability space) so that
 \beqnn
\<N(f)\>_t=\int_0^t\<\Phi(s)f,dB_s\>_{L^2([0,1]^2)},
 \eeqnn
where $\Phi(s)$ is linear map from $L^2(\IR)$ to $L^2([0,1]^2)$ so
that for $f\in L^2(\IR)$,
 \beqnn
(\Phi(s)f)(a,b)=\sqrt{\gamma(a,b)}\int_\mathbb{R}1_{\{a\leq
Y_s(y)\leq b\}}f(y)dy.
 \eeqnn
Let $\{h_j\}$ be a complete orthonormal system of the Hilbert space
$L^2([0,1]^2)$ and define random measure
$\{W(ds,da,db):s\ge0,~a,b\in[0,1]\}$ as
 \beqnn W([0, t]\times
A)=\sum_{j=1}^\infty\<1_A,h_j\>B_t^{h_j},\qquad t\ge0,~A\in
\mathscr{B}([0,1]^2).
 \eeqnn
It is easy to show that $\{W(ds,da,db):s\ge0,~a,b\in[0,1]\}$ is a
Gaussian white noise with intensity $dsdadb$. Furthermore,
 \beqnn
N_t(f)=\int_0^t\int_0^1\int_0^1\Big[\int_\mathbb{R}1_{\{a\leq
Y_s(y)\leq b\}}f(y)dy\Big]\sqrt{\gamma(a,b)} W(ds,da,db),
 \eeqnn
which implies the desired result.
 \qed

\subsection{Properties}

 \btheorem\label{t4.1}
Suppose that $\{X_t:t\ge0\}$ is the interacting Fleming-Viot process
with distribution-function process $\{Y_t:t\ge0\}$ solving
\eqref{4.1}. Then for any $f\in C^2(\mbb{R})$
 \beqlb\label{4.2}
X_t(f)
 \ar=\ar
X_0(f)+\frac{1}{2}\int_0^tX_s(f'')ds \cr
 \ar\ar
+\int_0^t\int_0^1\int_0^1[f(Y_s^{-1}(b))-f(Y_s^{-1}(a))]\sqrt{\gamma(a,b)}W(ds,da,db),
 \eeqlb
and for any $f\in B(\mbb{R})$, $\mbf{P}$-a.s.,
 \beqlb\label{4.3}
X_t(f) = X_0(P_tf)
+\int_0^t\int_0^1\int_0^1[P_{t-s}f(Y_s^{-1}(b))-P_{t-s}f(Y_s^{-1}(a))]
\sqrt{\gamma(a,b)}W(ds,da,db),
 \eeqlb
where $Y_t^{-1}(u)$ is defined by \eqref{3.11}. Moreover, for each
$t>0$ the random measure $X_t(dx)$ is absolutely continuous with
respect to the Lebesgue measure and the density $X_t(x)$ satisfies
 \beqlb\label{4.4}
 \ar\ar
X_t(x) =
\int_{\mbb{R}}p_t(x-y)X_0(dy)+\int_0^t\int_0^1\int_0^1\big[p_{t-s}(x-Y_s^{-1}(b))
\cr
 \ar\ar\qqquad\qqquad\qquad
-p_{t-s}(x-Y_s^{-1}(a))\big]
\sqrt{\gamma(a,b)}W(ds,da,db) \quad\mbf{P}\mbox{-a.s.}
 \eeqlb
 \etheorem
 \proof
The proof is similar to that of Theorem \ref{t3.2}.  By the Step 1 of the proof of Lemma
\ref{t4.2}, we obtain \eqref{4.2}. By a similar argument as
in the proof of \eqref{3.13}, we obtain \eqref{4.3} from \eqref{4.2}
immediately. By conditioning we may assume that $X_0$ is
deterministic and $X_0=\mu\in M_1(\mbb{R})$. It follows from
\eqref{4.3} that $\mbf{E}\{X_t(f)\}=\mu(P_tf)$ for $f\in
B(\mbb{R})$. Thus for $f\in B(\mbb{R})^+$,
 \beqnn
 \ar\ar
\mbf{E}\Big\{\int_{\mbb{R}}f(x)dx\int_0^tds\int_0^1da\int_0^1
[p_{t-s}(x-Y_s^{-1}(b))-p_{t-s}(x-Y_s^{-1}(a))]^2
\gamma(a,b)db\Big\} \cr
 \ar\ar\quad\le
2(2\pi)^{-\frac12}
\|\gamma\|\mbf{E}\Big\{\int_{\mbb{R}}f(x)dx\int_0^t(t-s)^{-\frac12}ds\int_0^1
p_{t-s}(x-Y_s^{-1}(a)) da\Big\} \cr
 \ar\ar\quad=
2(2\pi)^{-\frac12}
\|\gamma\|\mbf{E}\Big\{\int_{\mbb{R}}f(x)dx\int_0^t(t-s)^{-\frac12}ds\int_{\mbb{R}}
p_{t-s}(x-y) X_s(dy)\Big\} \cr
 \ar\ar\quad\le
\|\gamma\|\int_0^t(t-s)^{-\frac12}\mbf{E}\{X_s(P_{t-s}f)\}ds
=\|\gamma\|\mu(P_tf)\int_0^t(t-s)^{-\frac12}ds<\infty.
 \eeqnn
Now using \cite[Theorem 7.24]{Li11} we obtain
 \beqnn
X_t(f)=\int_{\mbb{R}}f(x)\check{I}_t(x)dx,
 \eeqnn
where $\check{I}_t(x)$ equals to the right-hand side of \eqref{4.4}.
This gives the desired result. \qed


\bigskip

\noindent{\small Li Wang}

\noindent{\small School of Sciences, Beijing University of Chemical
Technology, Beijing, P.R. China.}
\smallskip

\noindent{\small Xu Yang}

\noindent{\small School of Mathematics and Information Science,
Beifang University of Nationalities, Yinchuan, P.R. China.}
\smallskip

\noindent{\small Xiaowen Zhou}

\noindent{\small Department of Mathematics and Statistics, Concordia
University,1455 de Maisonneuve Blvd. West, Montreal, Quebec, H3G
1M8, Canada.}
\smallskip

\noindent{\small E-mails: {\tt wangli@mail.buct.edu.cn,
xuyang@mail.bnu.edu.cn, xiaowen.zhou@concordia.ca.}}


\noindent
\begin{thebibliography}{99}

\bibitem{BMP10}
Burdzy, K., Mueller, C. and Perkins, E. A. (2010): Nonuniqueness for
nonnegative solutions of parabolic stochastic partial differential
equations. \textit{Illinois J. Math.} \textbf{54}, 1481--1507.

\bibitem{ChenY15}
Chen, Yu-Ting (2015): Pathwise nonuniqueness for the SPDEs of some
super-Brownian motions with immigration. \textit{Ann. Probab.}
\textbf{43}, 3359--3467.

\bibitem{DLi12}
Dawson, D. A. and Li, Z. (2012): Stochastic equations, flows and
measure-valued processes. \textit{Ann. Probab.} \textbf{40},
813--857.



\bibitem{DD03} Delmas, J. F. and Dhersin, J. S. (2003): Super-Brownian motion with interactions. \emph{Probab. Theory Related Fields.}
\textbf{107}, 301--325.

\bibitem{EM90} El Karoui, N. and M\'{e}l\'{e}ard, S. (1990): Martingale measures
and stochastic calculus. \emph{Probab. Theory Related Fields.}
\textbf{84}, 83--101.


\bibitem{FV79} Fleming, W. H. and Viot, M. (1979): Some measure-valued markov processes
in population genetics theory.
 \emph{Indiana Univ. Math. J.} \textbf{28}, 817--843.

\bibitem{HLY14} He, H., Li, Z. and Yang, X. (2014):
Stochastic equations of super-L\'evy processes with general
branching mechanism. \textit{Stochastic Process. Appl.}
\textbf{124}, 1519--1565.

\bibitem{KX95} Kallianpur, G. and Xiong, J. (1995): \textit{Stochastic
Differential Equations in Infinite-Dimensional Spaces.} An extended
version of the lectures delivered by Gopinath Kallianpur as part of
the 1993 Barrett Lectures at the University of Tennessee, Knoxville,
TN.

\bibitem{KS88} Konno, N. and Shiga, T. (1998): Stochastic partial differential
equations for some measure-valued diffusions.  \emph{Probab. Theory
Related Fields.} \textbf{79}, 201--225.



\bibitem{Li11} Li, Z. (2011): \emph{Measure-valued Branching Markov Processes}.
Springer, Heidelberg.



\bibitem{MR93} M\'{e}l\'{e}ard, M. and Roelly, S. (1993): Interacting measure
branching processes: some bounds for the support. \emph{Stoch.
Stoch. Reports.} \textbf{44}, 103--121.



\bibitem{M85}
Mitoma, I. (1985): An $\infty$-dimensional inhomogeneous Langevin
equaion. \textit{J. Funct. Anal}. \textbf{61}, 342--359.

\bibitem{MMyP2014}
Mueller, C., Mytnik, L. and Perkins, E. (2014): Nonuniqueness for a
parabolic SPDE with $\frac34-\varepsilon$-H\"older diffusion
coefficients. \textit{Ann. Probab.} \textbf{42}, 2032--2112.

\bibitem{M02} Mytnik L. (2002): Stochastic partial differential
equation driven by stable noise. \textit{Probab. Theory Related
Fields.} \textbf{123}, 157--201.

\bibitem{MN15} Mytnik, L. and Neuman, E. (2015):
Pathwise uniqueness for the stochastic heat equation with H\"older
continuous drift and noise coefficients. \textit{Stochastic Process.
Appl.} \textbf{125}, 3355--3372.


\bibitem{MyP11} Mytnik, L. and Perkins, E. (2011): Pathwise uniqueness
for stochastic heat equations with H\"older continuous coefficients:
the white noise case. \textit{Probab. Theory Related Fields.}
\textbf{149}, 1--96.

\bibitem{MPS06} Mytnik, L., Perkins, E. and Sturm, A. (2006):
On pathwise uniqueness for stochastic heat equations with
non-Lipschitz coefficients. \textit{Ann. Probab.}, \textbf{34},
1910--1959.

\bibitem{MX14} Mytnik, L. and Xiong, J. (2015):
Well-posedness of the martingale problem for superprocess with
interaction. \textit{Illino. J. Math.} \textbf{59}, 485--497.

\bibitem{N14} Neuman, E. (2014):
Pathwise uniqueness of the stochastic heat equations with spatially
inhomogeneous white noise. Submitted. ArXiv:1403.4491.




\bibitem{Rei89} Reimers, M. (1989): One dimensional stochastic
    differential equations and the branching measure diffusion.
    \textit{Probab. Theory Related Fields.} \textbf{81}, 319--340.




\bibitem{RS13}
Rippl, T. and Sturm, A. (2013): New results on pathwise uniqueness
for the heat equation with colored noise. \textit{Electron. J.
Probab.} \textbf{18}, 1--46.




\bibitem{X13} Xiong, J. (2013): Super-Brownian motion as the unique
strong solution to an SPDE.  \emph{Ann. Proboa.} \textbf{41},
1030--1054.


\bibitem{Xiong13book}
Xiong, J. (2013): \textit{Three Classes of Nonlinear Stochastic
Partial Differential Equations}. World Scientific, Singapore.

\bibitem{XiongY}
Xiong, J. and Yang, X. (2016): Superprocesses with interaction and
immigration. \textit{Stochastic Process. Appl.} To appear.
\end{thebibliography}
\end{document}